\documentclass[12pt, reqno]{amsart}
\usepackage{ifthen,amsfonts,amsmath,amssymb,epic,eepic}
\usepackage{epsfig}

\makeatletter
\@addtoreset{equation}{section}
\makeatother

\theoremstyle{definition}

\def\fnum{equation}
\newtheorem{Thm}[\fnum]{Theorem}
\newtheorem{Cor}[\fnum]{Corollary}

\newtheorem{Lem}[\fnum]{Lemma}

\newtheorem{Pro}[\fnum]{Proposition}

\numberwithin{equation}{section}

\newcommand{\nn}{{\bf{n}}}

\newcommand{\dist}{{\text {dist}}}

\newcommand{\Hess}{{\text {Hess}}}
\def\RR{{\bf  R}}
\def\ZZ{{\bf  Z}}
\def\SS{{\bf  S}}

\newcommand{\lf}{{\ell_+}}

\newcommand{\Length}{{\text {Length}}}

\newcommand{\cF}{{\mathcal{F}}}

\newcommand{\cL}{{\mathcal{L}}}

\newcommand{\eqr}[1]{(\ref{#1})}

\begin{document}

\title[Geodesic laminations with closed ends on surfaces and Morse index]
{Geodesic laminations with closed ends on surfaces and Morse
  index; Kupka-Smale metrics}

\author{Tobias H. Colding}%
\address{Courant Institute of Mathematical Sciences and Princeton University\\
251 Mercer Street\\ New York, NY 10012 and Fine Hall, Washington
Rd., Princeton, NJ 08544-1000}
\author{Nancy Hingston}%
\address{Department of Mathematics\\
The College of New Jersey\\
Ewing, NJ 08628.}

\thanks{The first author was partially supported by NSF Grants DMS
  9803253 and DMS 0104453.}

\email{colding@cims.nyu.edu and hingston@TCNJ.EDU}

\maketitle

\section{Introduction}
Let $M^2$ be a closed orientable surface with curvature $K$
and $\gamma\subset M$
a closed geodesic.  The {\it Morse index} of
$\gamma$ is the index of the
critical point $\gamma$ for the length functional on the space of
closed curves, i.e., the number of
negative
eigenvalues (counted with multiplicity) of the second
derivative of length.  Since the second derivative of length at $\gamma$
in the direction of a normal variation $u\,\nn$ is
$-\int_{\gamma}u\,L_{\gamma}\,u$ where $L_{\gamma} \,u= u'' + K\,u$,
the Morse index is the number of
negative eigenvalues of $L_{\gamma}$.
(By convention,  an
eigenfunction $\phi$ with eigenvalue
$\lambda$ of $L_{\gamma}$ is a solution of
$L_{\gamma}\,\phi+\lambda\, \phi=0$.)  Note that if $\lambda=0$,
then $\phi$ (or $\phi\,\nn$) is a (normal) Jacobi field.
$\gamma$ is {\it stable}
if the index is zero.
The {\it index} of a noncompact geodesic is the dimension
of a maximal vector space of compactly supported variations for which the
second derivative of length is negative definite.  We also say that
such a geodesic
is {\it stable} if the index is $0$.

\vskip2mm
We give in this paper bounds for the Morse indices of a large
class of simple geodesics on a surface with a generic metric.  
To our knowledge these
bounds are the first that use only the generic hypothesis on the metric.  

\begin{Thm} \label{t:geolama}
For a generic metric on a closed surface, $M^2$, any geodesic
lamination with closed ends has
finitely many leaves and each leaf has finite Morse index.
\end{Thm}

Our second result is:

\begin{Thm} \label{t:geolam}
For a generic metric on a closed surface, $M^2$, there is a
bound for the Morse index of any collection of simple closed geodesics
for which each limit is a geodesic lamination with closed ends.
\end{Thm}

A {\it lamination}
on a surface $M^2$ is a collection $\cL$ of
smooth disjoint curves (called leaves)
such that
$\cup_{\ell \in \cL} \ell$ is closed.
Moreover, for each $x\in M$ there exists an
open neighborhood $U$ of $x$ and a coordinate chart, $(U,\Phi)$, with
$\Phi (U)\subset \RR^2$
so that in these coordinates the leaves in $\cL$
pass through $\Phi (U)$ in slices of the
  form $(\RR\times \{ t\})\cap \Phi(U)$.  
A {\it foliation} is a lamination for which 
the union of the leaves is all of $M$
and a {\it geodesic lamination} is a lamination whose leaves are geodesics.

If $\ell \in \cL$ is noncompact, then we set
\begin{equation} \label{e:ols1}
     \lf  = \cap_{s > 0} \overline{ \ell(s,\infty) } \, .
\end{equation}
Since $\lf$ is the intersection of nonempty nested closed sets it is
closed and nonempty since $M$ is
compact.    Since $\cup_{\ell \in \cL}
\ell$ is closed, $\lf\subset \cup_{\ell \in \cL}
\ell$.  Likewise we define $\ell_-$.
A leaf $\ell\in \cL$ is said to be isolated if for some $x\in \ell$
(hence all $x\in \ell$) there exists $\epsilon=\epsilon (x)>0$ such that
$B_{\epsilon}(x)\cap \tilde\ell=\emptyset$ for all
$\tilde\ell\in \cL\setminus \{\ell\}$.
Note that $\ell_-$, $\ell_+$
consist of nonisolated leaves.  
We say that a geodesic lamination $\cL$ has {\it closed ends} if for each
noncompact leaf $\ell\in \cL$ both $\ell_+$ and $\ell_-$ are closed
geodesics.

We will equip the space of metrics on a given manifold with the
$C^{\infty}$-topology.  A subset of the set of metrics on a
given manifold is said
to be {\it residual} if it is a countable intersection of open dense
subsets.
A statement is said to hold {\it for a generic metric} if it holds for
all metrics in a residual set.  
The conclusions of Theorems \ref{t:geolama} and \ref{t:geolam} 
are true for a residual set of metrics
that we call {\it Kupka-Smale} (KS-metrics).  
This hypothesis on the metric has a 
natural 
interpretation in both the dynamical systems and the variational contexts.  
Here are two versions of our hypothesis:

\vskip2mm
\underline{KS-metric (dynamical version)}:\\
(1) Every simple closed orbit of the geodesic 
flow whose Poincar\'e map has real eigenvalues is hyperbolic.\\
(2) Every intersection of stable- and unstable-manifolds at a 
simple geodesic is transverse.

\vskip2mm
\underline{KS-metric (variational version)}: 
Let $\gamma$ be a simple geodesic.\\
(1) If $\gamma$ is periodic, there is no periodic Jacobi field without  
zeroes.\\
(2) If $\gamma$ is noncompact and has closed ends, there is no bounded 
Jacobi field without zeroes. 

\vskip2mm
 The above two conditions are equivalent under the additional 
condition that the metric is bumpy.  
A metric on a surface is {\it bumpy} if each  closed geodesic
is a nondegenerate critical point, i.e.,
$L_{\gamma} u = 0$ implies $u\equiv 0$.  Bumpy metrics are generic,
\cite{Ab}, \cite{An}. For convenience we will prove the conclusions of
 Theorems 
\ref{t:geolama} and \ref{t:geolam} for these ``bumpy KS-metrics".  For 
further discussion of these metrics, see Section \ref{s:s3} where we
 show that the set of bumpy KS metrics contain a residual set.  

 Here is the idea of the proof of Theorem 
\ref{t:geolama}:    
Unstable closed leaves in a fixed lamination are always isolated; 
thus we need to show that noncompact leaves are isolated and have 
finite index.  In Theorem \ref{t:geolama} part 1 of the 
KS-condition is applied to get a nice structure on the
 ends of noncompact leaves, and to ensure that these ends are isolated and
 countable.  If an end leaf is closed and hyperbolic, this  structure is
 striking:  On each side of the end (limit) leaf there are two smooth circles
 of geodesics, each spiraling toward the limit leaf, one in each direction
 on each side; see fig. 1 and Corollary \ref{c:stablemfldf}. 
Each circle foliates a tubular neighborhood of the given side
 of the end.  These circles of noncompact geodesics are the stable- and
 unstable-manifolds of the end leaf, when viewed as a closed orbit of the
 geodesic flow on the unit tangent bundle of the surface.  The second part of
 the KS condition ensures that, in a given lamination, leaves limiting on a
 given pair of ends are isolated, as the corresponding  circles intersect
 transversely in a local (two-dimensional) section of the flow; see 
fig. 2. These
 noncompact leaves have finite index since index ``stops accumulating" once
 they get close to the stable, hyperbolic, ends.

 To prove Theorem 
\ref{t:geolam} we first extract a converging subsequence of the given
sequence of simple closed geodesics.  
The limit is easily seen to be a geodesic
lamination with multiplicities; 
see the discussion preceding Proposition \ref{p:finitep2}. 
- On long stretches the geodesics in the subsequence
will mimic the behavior of the limit lamination.  
By assumption this limit lamination has 
closed ends.  
The transversality of the intersection of the stable- and
unstable-manifolds will then allow us to conclude finiteness of the
indices for 
the converging subsequence of geodesics. 

\medskip
\centerline{\epsfig{figure=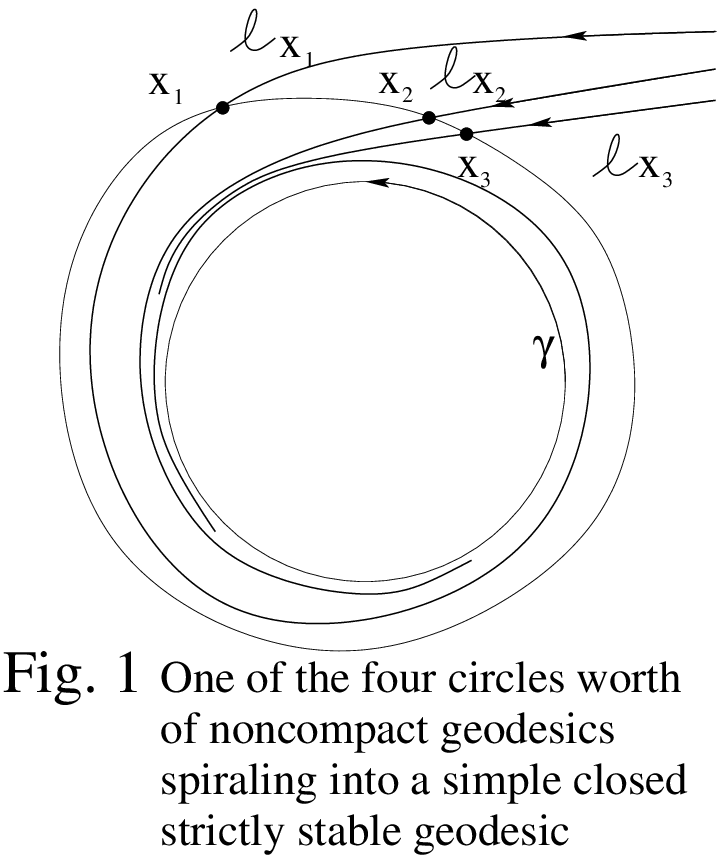}\hspace{0.5cm}
\epsfig{figure=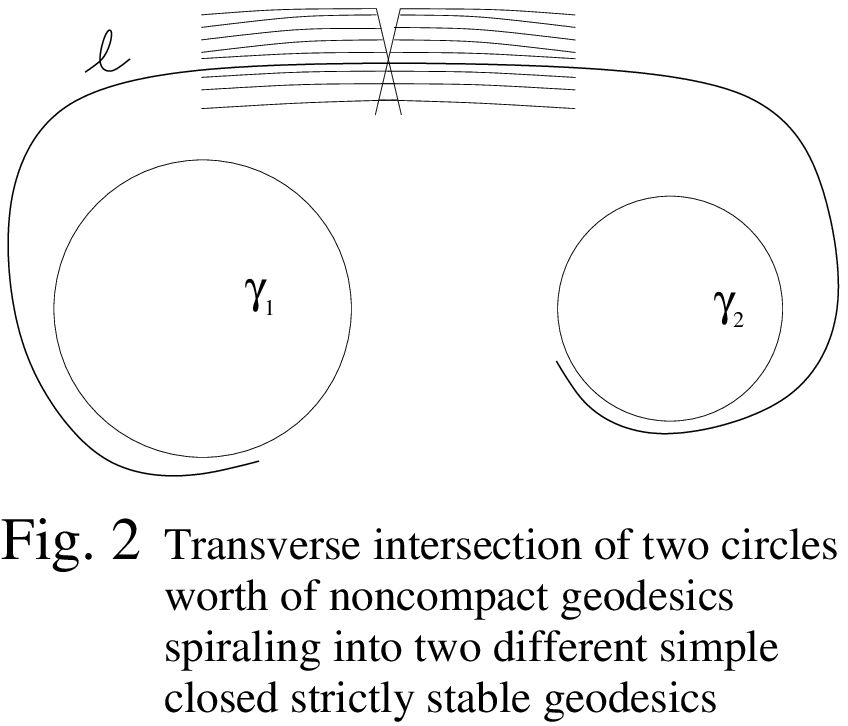}}
\medskip

It is easy to see that the first part of Theorem  \ref{t:geolama} is
false without the word ``generic''.  One can construct a surface of
revolution that has 
a geodesic lamination with infinitely many leaves.  In this example all
the leaves are stable and have closed ends.  However in 
\cite{CH1} we showed that on any $M^2$,
there exists a metric with a geodesic lamination
with closed ends and infinitely many unstable leaves.  
Moreover, there exists such a metric
which has no bound for the
index of all simple closed geodesics.

General geodesic laminations on surfaces need not have closed ends; 
consider for instance a flat square torus with
the foliation consisting of lines with a common irrational slope.
In fact, on any surface there are (bumpy) metrics and geodesic laminations
without closed ends:

\begin{Thm}  \label{t:notclosed}
On any surface $M^2$, there exists an open 
(nonempty) set of metrics having   
geodesic laminations without closed ends.  These laminations are
limits of sequences of simple closed geodesics.
\end{Thm}

Our interest in whether the Morse index is bounded for simple closed 
geodesics on surfaces comes in part 
from its connection with the spherical space form problem; 
see \cite{PiRu}, \cite{CM2} where Pitts and Rubinstein ask for such a 
bound for embedded minimal tori for a sufficiently large class of 
metrics on $\SS^3$.  Clearly obstructions to Morse index bounds for simple 
closed geodesics on surfaces give obstructions to Morse index bounds for 
embedded minimal tori on $3$-manifolds 
(the most immediate generalization of simple
closed geodesics on surfaces to $3$-manifolds is embedded minimal tori 
with uniform curvature bounds).  In addition 
new obstructions occur; \cite{HaNoRu}, \cite{CD}.  
We believe that many of the ideas of this 
paper can be used to give bounds for the Morse indices of 
geodesics on surfaces and 
on embedded minimal tori (or more generally fixed genus) with uniform
curvature bounds on $3$-manifolds with generic metrics.  In fact, 
the arguments given here should be useful even for surfaces 
without curvature bounds; for instance on 
closed $3$-manifolds with positive scalar curvature any complete 
stable minimal surface is necessarily 
compact and in fact either topologically $\SS^2$ or $\bf{RP}^2$.  Thus 
``ends'' of embedded minimal annuli (even without curvature bounds) 
in such manifolds are closed; see 
\cite{CD}, \cite{CM3} for more discussion on this.

\vskip2mm
Throughout this paper 
$M^2$ is a closed orientable surface with a Riemannian metric, $\cL$
is a geodesic lamination,
and when $x\in M$, $r_0>0$, and
$D\subset M$, then we let $B_{r_0}(x)$ denote the ball of
radius $r_0$ centered at $x$ and $T_{r_0}(D)$
the $r_0$-tubular neighborhood of $D$.  Moreover, if $x,$ $y\in M$, then
$\gamma_{x,y}:[0,\dist_M(x,y)]\to M$ will denote a minimal geodesic
from $x$ to $y$.  Whenever we look at a single geodesic it will always
be assumed to be parameterized by arclength.

\vskip2mm
We are grateful to Camillo De Lellis for making the illustrations.

\section{Geodesic laminations on surfaces}  \label{s:s1}

We will often implicitly use the following simple fact:  If
$\gamma\subset M^2$ is a simple closed geodesic, then there exists
$\delta=\delta (\gamma)>0$ such that the nearest point projection
$\Pi_{\gamma} : T_{\delta} (\gamma) \to \gamma$ is well defined.
Moreover, if
$\tilde \gamma :[0,1]\to T_{\delta}(\gamma)$ is a geodesic,
then
\begin{equation}  \label{e:pigamma}
||d \Pi_{\gamma} |_{\tilde \gamma}| - 1 | < \psi (\delta)
\end{equation}
where $\lim_{\delta \to 0} \psi (\delta)=0$.  Note that this
just says that the geodesics $\gamma$ and $\tilde \gamma$ are nearly
parallel.
If $\gamma$ is oriented, then we say that
$\tilde \gamma:[0,1]\to T_{\delta}(\gamma)$ has
the same orientation as $\gamma$ if
$|d\Pi_{\gamma}\tilde \gamma'-\gamma'|<\psi (\delta)$.
We will most often assume that this is the case.

We will assume in what follows some knowledge of Jacobi fields and indices
of geodesics; see e.g. \cite{Kl} or \cite{Sp}.  Three facts will be
particularly important:\\
1) A Jacobi field $J$ is uniquely determined by the
values $(J(t),J'(t))$ for any $t$.\\
2) A complete geodesic (closed or
noncompact) is stable if and only if it has no Jacobi field with more than
one zero. (In the closed case we of course mean Jacobi field with the same
periodicity as the geodesic.)\\
3) A noncompact geodesic has finite index if and only if
there is a bound for the number of zeros of any nontrivial
Jacobi field along it.

To show that certain geodesics are stable (or have bounded Morse index) it
is sometimes useful to apply the following standard fact:  A
Schr\"odinger operator $L\,u=u''+K\,u$ is nonpositive
($-L\geq 0$) if it has a positive supersolution $\phi$
(that is $\phi>0$ and $L\,\phi\leq 0$).  This follows since
$-(\log \phi )''\geq K+|(\log \phi)'|^2$ and hence if $f$ is a
compactly supported function, then
integration by parts and the Cauchy-Schwarz inequality yields
\begin{align} \label{e:poscri}
\int f^2\,K+\int f^2\,|(\log \phi )'|^2
\leq & -\int f^2\,(\log \phi)''\\
= & 2 \int f\,f'(\log \phi)'
\leq \int f^2\,|(\log \phi)'|^2+\int (f')^2\, .\notag
\end{align}
Thus, $-\int f\,L\,f=-\int f\,(f''+K\,f)\geq 0$.

\begin{Lem}  \label{l:stablengbh}
Let $\gamma$ be a strictly stable ($-L_{\gamma}>0$) 
simple closed geodesic on $M^2$.
There exists $\delta=\delta (\gamma)>0$, such that any geodesic
segment contained in
$T_{\delta}(\gamma)$ with length $\geq 1$ is stable.
\end{Lem}

\begin{proof}
Since $\gamma$ is strictly stable, then $-L_{\gamma}> 0$ so if we let
$\lambda_1$ be the first eigenvalue and $\phi$ a corresponding
eigenfunction, then $\lambda_1>0$ and $\phi^2>0$.  In particular
$-L_{\gamma}|\phi|=\lambda_1\,|\phi|>0$.  Let
$\tilde \gamma\subset T_{\delta}(\gamma)$ be a geodesic segment with length
$\geq 1$ and
set $\tilde \phi=\phi\circ \Pi_{\gamma}$, then $|\tilde \phi|>0$ and
(by \eqr{e:pigamma})
$-L_{\tilde \gamma}|\tilde \phi|>0$.  The lemma now follows from the
remarks preceding it.
\end{proof}

Let $\gamma$ be a closed geodesic with universal cover $\tilde
\gamma$.  Then\\
1). $\gamma$ is stable if and only if $\tilde \gamma$
has no Jacobi field $J:\RR\to \RR$ with $2$ zeroes.\\
2). $\gamma$ is strictly stable if and only if $\gamma$ is stable and 
$\tilde \gamma$ has no periodic Jacobi field.\\
The ``only if'' part of each statement follows by contradiction from
the simple argument of Lemma \ref{l:stablengbh} when applied to a
first eigenfunction $\phi$, with $f$ the restriction of $J$ for 1). 
to an interval
between two zeroes, and for 2) to one period of $J$.  
The ``if'' part of 2) is clear; to see the ``if'' part of 1). we argue
using \eqr{e:poscri}:  Let $\chi_n$ be the cutoff function 
(we may assume that $\text{Length}(\gamma)=1$) 
\begin{equation}
   \chi_n (t)=
\begin{cases}
    1 &\hbox{ for } |t|\leq n^2\, ,    \\
    1-(|t|-n^2)/n &\hbox{ for } n^2<|t|\leq n^2+n\, ,\\
    0 &\hbox{ otherwise }  \, .\\
\end{cases}
\end{equation}
If $f$ is a
function on $\gamma$, $\tilde f$ its lift to
$\tilde \gamma$, set $f_n=\tilde f\, \chi_n$, and $\phi_n=J_n$
(where $J_n$ is a Jacobi field with
$J_n(-n^2-n)=0$ and $J_n|(-n^2-n,\infty)>0$.  
Then by \eqr{e:poscri} for $f_n$
\begin{align}  \label{e:converse}
2\,n^2\int_{\gamma} K\,\tilde f^2&\leq 
\int K\,\tilde f^2\,\chi_n^2\leq \int [(\tilde f\,\chi_n)']^2\notag\\
&\leq \int [(\tilde f)']^2\,\chi_n^2
+\int_{n<|t|\leq n^2+n}[(\tilde f)^2+(\chi_n ')^2]
+\int (\tilde f)^2\,(\chi_n')^2\\
&\leq 2\,n\,(n+2)\int_{\gamma}(f')^2
+\int_{n<|t|\leq n^2+n}(\chi_n ')^2
+\int (\tilde f)^2\,(\chi_n')^2\, .\notag
\end{align}
Dividing \eqr{e:converse} through by $2\,n^2$ and 
letting $n\to \infty$ gives 
$-\int_{\gamma} f\,L\,f=-\int_{\gamma} f\,(f''+K\,f)\geq 0$.

We will also use the following two well known facts;
see e.g. \cite{Kl} or \cite{Sp}.  One can also see 4) below
by arguing as in 1) and 2) above:\\
3) A Jacobi field $J$ is uniquely determined by the
values $(J(t),J'(t))$ for any $t$.\\
4) A noncompact geodesic is stable if and only if it has no
Jacobi field with more than
one zero and it has finite index if and only if
there is a bound for the number of zeros of any nontrivial
Jacobi field along it.

Note that 1) and 4) together imply that $\gamma$ is stable if
and only if $\tilde \gamma$ is.

\begin{Lem}  \label{l:nonisolated}
Any nonisolated leaf $\ell$ of $\cL$ is stable.
\end{Lem}

\begin{proof}
If $\ell$ had a Jacobi field with $2$ zeroes, then so
would every sufficiently nearby (in the unit tangent
bundle) geodesic.  But between $\ell$ and any nearby
geodesic $\ell_{i}$ which does not intersect $\ell$, we
can find (using $\ell$ and $\ell_{i}$ as barriers) a
stable geodesic $\eta$ which has no Jacobi field with $2$ zeroes.
\end{proof}

    From the definition of a lamination and Lemma \ref{l:nonisolated}
one easily shows:

\begin{Lem} \label{l:ls1}
Each $\lf$ is connected (as a subset of $M$).
Moreover, given $x \in \lf$, then $x \in \ell_0 \subset \lf$
for some $\ell_0 \in \cL$;
$\ell_0$ is said to be a {\it limit
leaf} and is stable.
\end{Lem}

We will need the following:

\begin{Lem} \label{l:+closed}
If $\ell\in \cL$ is
noncompact and $\ell_+$ contains a closed geodesic $\sigma$,
then $\ell_+=\sigma$ and $\ell$ spirals
monotonically toward $\ell_+$.
\end{Lem}

\begin{proof}
Observe first that once $\ell$ get into a small tubular neighborhood of
$\sigma$, then one of the two ``directions'' of $\ell$ must be completely
contained in a small tubular neighborhood of $\sigma$.  Now using that
$\sigma\subset \ell_+$ it follows from this that the ``forward direction''
of $\ell$ is actually contained in a small tubular neighborhood of
$\sigma$ and hence (again since $\sigma\subset \ell_+$)
must spiral towards $\sigma$ monotonically.
\end{proof}

We say that a geodesic $\gamma_2:[0,k_2]\to M^2 $ can be written as a {\em
normal graph} over a geodesic $\gamma_1 :[0,k_1]\to M$
by a function $u$ (on $[0,k_1]$) if there
is a diffeomorphism $\alpha :[0,k_1]\to [0,k_2]$ such that for all
$t\in [0,k_1]$, then
\begin{equation}
\gamma_2 (\alpha (t))=\exp_{\gamma_1(t)}(u(t)\,\nn_{\gamma_1}(t))\, .
\end{equation}

The next corollary follows from Lemma \ref{l:+closed}.

\begin{Cor} \label{c:s2graph}
Suppose that $\cL$ has closed ends.  
If $\ell_i\,,\ell\in \cL$ and $\ell$ is noncompact
with $\ell_i'(0)\to \ell' (0)$, then for $i$ sufficiently large
$\ell_{i_{\pm}}=\ell_{\pm}$ and $\ell_i$ is a normal graph over $\ell$.
\end{Cor}

Recall that we equip the space of $C^{\infty}$ metrics on a closed surface
$M^2$ with the
$C^{\infty}$-topology and we write $g_i\to g$ if $|g-g_i|_{C^{\infty}}\to
0$.
Most of the next lemma will be needed only in Section \ref{s:s3}.

\begin{Lem}  \label{l:vary}
(Lemma B.1 of \cite{CH1}). 
Suppose that the metric $g$ on $M$ is bumpy.  For each $L>0$,
there exists at most finitely many closed geodesics of length $<L$.
Moreover, if $L$ is not equal to the
length of any closed geodesic in $g$, then in a
neighborhood of $g$ each metric has precisely as many (simple)
closed stable geodesics of length $<L$ as $g$.  Finally, if
$g_i\to g$ and $\{\gamma_{i,k}\}$, $\{\gamma_k\}$
are the (simple)
closed stable geodesics in
$g_i$, $g$, respectively, of length $<L$,
then $\gamma_{i,k}\to \gamma_k$ for $i\to \infty$ and
each $k$.
\end{Lem}

     In the remainder of this section the metric on $M^2$ is bumpy
metric and $\cL$ has closed ends.  By Lemma \ref{l:vary}, 
$\cL$ contains finitely many closed leaves $\eta_1 , \dots, \eta_m$; these
are the only limit leaves.  Choose  $\epsilon > 0$ so that:
\begin{equation}  \label{e:eprest}
    \dist (\eta_j , \eta_k) > 2 \,
    \epsilon {\text{ for }} j \ne k \, ;
\end{equation}
for $1 \leq j \leq m$ and $\ell \in \cL$,  $T_{\epsilon} (\eta_j) \cap \ell$
is graphical over $\eta_j$; and
$\Pi : T_{\epsilon} (\cup \eta_j) \to \cup \eta_j$ is smooth.
Using the local product structure, there exist
$C > 0$ and
$S$ with $\partial S$ smooth so that
\begin{equation} \label{e:mys}
    \cup T_{\epsilon / C}(\eta_j) \subset S \subset
    \cup T_{\epsilon}(\eta_j)
\end{equation}
and $\partial S$ intersects $\cL$ transversely.
$S_j$ denotes
the component of $S$
containing $\eta_j$.

\begin{Cor} \label{c:short}
Suppose that the metric on $M$ is bumpy and $\cL$ has closed ends.  
Let $\eta_j$, and $S$ be as above.
There exists
$\rho < \infty$
so that, for each $\ell \in \cL$, each
    component $\alpha$ of
$\ell \setminus S $
has $\Length (\alpha) \leq \rho$.
\end{Cor}

\begin{proof}
This follows by compactness.
\end{proof}

Using the local product structure,
Corollary \ref{c:short} implies that there
are collections $T_j$, $1 \leq j \leq \bar{m}$,
of components $\alpha$ of
$\cup_{\cL} \ell \setminus S$
so that each $\alpha \in T_j$ spirals between
$\eta_{1_j}$ and $\eta_{2_j}$.  Hence,
\begin{equation}    \label{e:epsbar}
       \min_{j_1 \ne j_2} \, \dist (\cup_{\alpha \in T_{j_1}} \alpha ,
       \cup_{\alpha \in T_{j_2}} \alpha ) = \bar{\epsilon} > 0 \, .
\end{equation}

\section{Morse index bounds for bumpy Kupka--Smale metrics}
\label{s:s2}

In this section we discuss the stability of the leaves of a geodesic
lamination.  The first goal is to prove Corollary \ref{c:finitel},
which says that on a closed surface $M^2$ with a bumpy metric,  
a nonisolated leaf in a geodesic lamination
with closed ends 
implies a non-transverse intersection of two circles which will be described
below.  From this corollary we will then be able to give a
condition (KS-metrics) on a metric on $M$ which will imply that all geodesic 
laminations with
closed ends have
finitely many leaves, and that there is a bound for the Morse index of 
simple closed geodesics; see Propositions \ref{p:finitep}
and \ref{p:finitep2}.  In Section \ref{s:s3} we will see that this
condition is generically satisfied.

If $\gamma:[t_1,t_2]\to M^2$ is a geodesic,
then we let $P_{t_{2},t_{1}}$ denote
the (relative) linear Poincar\'e map which describes to first order
how nearby geodesics advance along $\gamma$.
That is, if $(a,b)\in \RR^2$, then
$P_{t_{2},t_{1}}(a,b)=(J(t_{2}),J'(t_{2}))$ where $J$ is
the
Jacobi field on $\gamma$ with $J(t_1)=a$ and $J'(t_1)=b$.
Note that $P_{t_2,t_1}=P_{t_{2},t}\,P_{t,t_1}$.
Set
\begin{gather*}
R_{\gamma}(t)=
\begin{pmatrix}
0 & 1 \\ -K(\gamma (t)) & 0
\end{pmatrix}
\, .
\end{gather*}
By the Jacobi equation $\frac{d}{dt}P_{t,t_1}=R_{\gamma}(t)\,P_{t,t_1}$;
since
$\text{Tr}\,(R_{\gamma})=0$ and $P_{t_1,t_1}$ is the identity
it follows that
$P_{t,t_1}\in \text{SL}(2,\RR)$.  Observe that if
$\gamma:[0,s_{\gamma}]\to M$ is closed, then
$P_{\gamma}=P_{s_{\gamma},0}$ is the usual linear Poincar\'e map.

The existence of a Jacobi field along $\gamma$ with zeroes at $t_1$ and
$t_2$ is equivalent to the fact that $P_{t_2,t_1}$ (as a linear
map from $\RR^2$ to itself) takes the $y$-axis to itself;
thus we will want to keep
an eye on the $y$-axis as $P_{t,t_1}$ acts on $\RR^2$.   Note that the $1$ in
the upper right corner of $R_{\gamma}$ means that if we watch the motion of
a vector $(a,b)$ under $P_{t,t_1}$, then at a time $t_2$ when
the vector hits the $y$-axis (i.e. when $P_{t_2,t_1}(a,b)
=(0,y)$), the vector is moving clockwise, i.e. $\frac{d}{dt}
P_{t,t_1}(a,b) =(y,0)$. (This agrees with common sense:
If $J(t_2)=0$, and $J'(t_2)=y$, then $J(t)$ has the same sign as $y$ for
$t>t_2$.)

    If the metric on $M$ is bumpy, and $\gamma$ is a closed limit leaf
of a geodesic lamination, then $\gamma$\\
1). is simple,\\
2). has no Jacobi field $J:\RR\to \RR$ with $2$ zeroes,\\
3). has no periodic Jacobi field.\\
If $\gamma$ is closed, it is clear by the above discussion and continuity
that if $P_{\gamma}=P_{s_{\gamma},0}$ does not have a positive real
eigenvalue (i.e. if $P_{\gamma}$ does not fix the direction of some vector
in $\RR^2$), then the path $P_{t,0},0\leq t\leq s_{\gamma}$, rotates each
vector in $\RR^2$ clockwise a positive amount.  However since
$P_{ms_{\gamma}+t,0}=P_{t,0}\,P_{s_{\gamma},0}^m$, in that case eventually
the $y$-axis will be mapped to itself, causing a Jacobi field with $2$
zeroes. Thus if $\gamma$ is simple closed and
strictly stable, the eigenvalues of $P_{\gamma}$
are of the form $\lambda$ and $1/\lambda$
where $0<\lambda<1$; it follows that $P_{\gamma}$
has a basis (not necessarily orthogonal)
of eigenvectors. (Note that $1$ cannot be an eigenvalue by 3).)

A local
    section $\Sigma$ of the geodesic flow along $\gamma$ is obtained as
follows; \cite{Bi}, \cite{MS}.
Pick $t_{0}$, and construct a geodesic $\tau$ on $M$
    transverse to $\gamma '(t_{0})$ at $\gamma (t_{0})$.
Let $\Pi$ be the projection from the unit tangent
    bundle $T_{1}M$ onto $M$.  The surface
    $\Sigma\subset T_{1}M$ is the intersection of a
    neighborhood of $\gamma '(t_{0})$ with the set
$\Pi ^{-1}(\tau)$.  Each point in $\Sigma$ corresponds (by
giving an initial
tangent vector) to a
    geodesic near $\gamma$.  If $\gamma$ is closed, by
following
    the geodesics around we obtain the ($C^1$)
Poincar\'{e} map
    ${\mathcal P}:\Sigma \rightarrow \Sigma$ with fixed point
$\gamma '(0)$. (Strictly speaking we will need to make
the domain of
${\mathcal P}$
smaller in order to get the range inside $\Sigma$.)
    The derivative of ${\mathcal P}$ at $\gamma '(0)$ is the
linear
    Poincar\'{e} map $P=P_{\gamma}$.

We will use without proof
the
following lemma, which says that an appropriate limit
of geodesics is a
Jacobi
field.  We decline to put a topology on the set of
geodesics on $M$.
However, very loosely speaking, if we think of the
space of Jacobi fields
along
$\gamma$ as the tangent space to the set of geodesics
at $\gamma$, then the
lemma says that a neighborhood of $\gamma$
in the set of geodesics is diffeomorphic to a
neighborhood of $\gamma
'(t_{0})$ in $\Sigma$. The ``diffeomorphism" takes a
geodesic $\sigma$ to
its
tangent
vector $\sigma '(t)$ at the time $t$ when it crosses
$\tau$, and a Jacobi
field along the
geodesic to the values $(J(t),J'(t))$ at that time.

\begin{Lem}  \label{l:linearna}
Let $\gamma$ be a geodesic, and for $i\geq 1$ let
$u_{i}(t)$ be the normal graph over $\gamma$ of a geodesic
$\gamma_{i}$. Assume that $\lim_{i\to\infty}\|(u_{i}(0),u_{i}'(0))\|=0$,
and that the limit
\begin{equation}
\lim_{i\rightarrow\infty} (u_{i}(0),u_{i}'(0))/\|(u_{i}(0),u_{i}'(0))\|
\end{equation}
exists. Fix $t_{0}$. Then $\lim_{i\rightarrow \infty}
u_{i}/\|(u_{i}(t_{0}),u_{i}'(t_{0})\|$ exists, and represents a Jacobi field
$J$ with
\begin{equation}  \label{e:linearna}
( J(t),J'(t)) \:=\: \lim_{i\rightarrow
\infty}( u_{i}(t),u_{i}'(t))/\|(u_{i}(t_{0}),u_{i}'(t_{0}))\|
\end{equation}
for all $t$.
Conversely, any Jacobi field is the limit
of a $1$-parameter family of
geodesics (though it may be that none of these
geodesics is the normal
graph of a function $u(t)$ defined for all $t$).
\end{Lem}

Suppose now that $\gamma$ is simple closed and strictly stable.  The
lemma that follows
says that the Poincar\'e map ${\mathcal P}$ has the same
behavior as its
derivative $P$:  It has one contracting direction
(eigenspace for
$\lambda$) and one expanding direction (eigenspace
for $1/\lambda$).

\begin{Lem} \label{l:stablemanif}
(See \cite{HiPu}). Let $\gamma$ be simple closed and strictly stable.  If
$\Sigma$ is sufficiently small, then
there is a $C^1$ curve $Y$ through the origin in
$\Sigma$ with the property that, for all $x\in\Sigma$,
\begin{equation}
\lim_{n\to \infty}{\mathcal P}^{n} x=0
\:\Leftrightarrow\: x\in Y\:\Leftrightarrow
\:\forall n>0\:{\mathcal P}^{n}x\in\Sigma \, .
\end{equation}
\end{Lem}

    Thus a geodesic $\tau$ near $\gamma$ has $\tau_{+}=\gamma$
if and only if the point $x$ in $\Sigma$
corresponding to $\tau$ lies on $Y$. $Y$ is called the
{\em stable manifold} of ${\mathcal P}$.

\begin{Cor} \label{c:}
Let $\gamma$  be simple closed and strictly stable.
If $\ell$ is a noncompact simple geodesic with $\ell_{+}
=\gamma$ (or $\ell=\gamma$), then $\ell$ has a unique Jacobi field $J_{+}$
with the normalization $X_{+}(0)
\in\{(\cos\theta ,\sin\theta )\mid -\pi /2 < \theta \leq
\pi /2\}$, where $X_{+}(t)=(J_{+}(t),J_{+}'(t))$, and
with $\| J_{+}(t)\|$ bounded for $t>0$.  For this
vector field there exists $C$ and $\alpha >0$ such that
for $t >0$
\begin{equation}
\|X_{+}(t)\| \leq C\exp (-\alpha t)\, .
\end{equation}
\end{Cor}

\begin{proof}
Let $x$ be a point in $\Sigma$ representing $\ell$;
this means $x=\ell '(t_{0}).$ A vector in the tangent
space to $\Sigma$ at $x$ represents (by giving the
initial values $J(t_{0})$ and $J'(t_{0})$) a Jacobi
field along $\ell$.
By Lemma \ref{l:linearna}, the derivative $d_{x}{\mathcal P}$
describes how the Jacobi field $J$ advances along one
loop of $\ell$. Using \eqr{e:linearna} it is clear that a
tangent
vector to the stable manifold $X$ at $x$ represents a
Jacobi field with the desired property.
\end{proof}

\begin{Cor} \label{c:stablemfldf}
(See fig. 1).  
Let $\gamma$  be simple closed and strictly stable.
Then there are four ``circles" of
noncompact
geodesics limiting on $\gamma$.  That is, on each side
of $\gamma$ in $M$,
and for each
orientation of $\gamma$ there is a $C^1$ map
$\SS^{1}\rightarrow T_{1}M$ which gives a bijection
between the circle $\SS^{1}$ and the set of geodesics
$\ell$ with $\ell_{+}=\gamma$ which limit on $\gamma$
from the given side
of $M$ with the given orientation.
\end{Cor}

\begin{proof}
Implicit in the statement of the corollary is the map
from $T_{1}M$ to the set (untopologized) of geodesics
on $M$.  One way to parameterize the circle is to use
the segment of the stable manifold $Y$  between two
consecutive points $x$ representing a single geodesic
$\ell$ with $\ell_{+}=\gamma$.  This does not quite
work, as the endpoints of the segment correspond to
different starting points on the same geodesic; however
it is clear how to reparameterize along the segment in
order to get the ends to match up. Once this is done,
the image of $\SS^{1}$
in $T_{1}M$ will be a circle close to the curve of
tangents of $\gamma$,
whose image in $M$ lies on the given side of
$\gamma$.
\end{proof}

If a geodesic $\ell$ lies in one of these four circles of
geodesics  given by Corollary \ref{c:stablemfldf},
then the vector field $J_{+}$ along $\ell$ can be
thought of as a tangent vector to the circle.

\begin{Cor} \label{c:foliation}
Let $\gamma$  be simple closed and strictly stable.
There exists a neighborhood $T$ of $\gamma$ such
that for all $x\in T\setminus \{\gamma\}$
there is a unique (maximal)
geodesic $\ell_x:(a,\infty)\to T$
with the same orientation
as $\gamma$ and $x\in \ell_x$.
Moreover, $a>-\infty$,
$\partial \ell_x\in \partial T$, $\ell_x\subset T\setminus \{\gamma\}$,
$\ell_x$ is simple, $(\ell_x)+=\gamma$, and
$\cF=\{\ell_x\}_{x\in \partial T}\cup \{\gamma\}$
is a geodesic foliation of $T$.
\end{Cor}

\begin{proof}
Let $T$ have as its boundary the image in $M$ of the
circle in $T_{1}M$
given in the previous corollary. (This will need to be
done once on each side
of $\gamma$.) To get the foliation structure, and the
simplicity of $\ell_{x}$, we use the fact that the
derivative of the
composite
\begin{equation}
Y\longrightarrow \Sigma \underset{\Pi}{\longrightarrow}M
\end{equation}
at the point
$\gamma '(t_{0})$ is nontrivial. The latter fact can be
seen as follows: A tangent vector to $Y$ at $\gamma '(t_
{0})$ gives the initial values $(J(t_{0}),J'(t_{0}))$ of
a Jacobi field which returns after one loop around
$\gamma$ as a multiple $1/\lambda$ of itself.  Since
$\gamma$ is stable and thus has no Jacobi field with
two zeroes, $J(t_{0})\neq 0$ and thus the image under
$\Pi$ is nonzero.
\end{proof}

There is a direct way (using the appendix of \cite{CH1}) of
getting Corollaries \ref{c:stablemfldf}, \ref{c:foliation} 
without appealing to Lemma \ref{l:stablemanif}.
Namely, by appendix A of \cite{CH1} there exists a strictly convex
function $F$ defined in a neighborhood $\{F\leq \epsilon\}$
(where $\epsilon>0$ is sufficiently small) of $\gamma$.  (In fig. 1
the curve circling $\gamma$ is meant to represent a level set of $F$.)  
Note that each
side of $\gamma$ in this neighborhood is convex and homeomorphic to a
cylinder.  A straightforward convergence argument show
that for each $x\in \{F=\epsilon\}$
there exists a simple stable geodesic
$\ell_x\subset \{F\leq \epsilon \}$ with $x\in \ell_x$ and
$(\ell_x)_+=\gamma$ as
in Corollary \ref{c:foliation}.  That $\ell_x$, $\ell_y$ do not
cross (and that $\cup_x\ell_x=\{F\leq \epsilon\}\setminus \{\gamma\}$)
follows easily from Lemma \ref{l:linearna} using the linear
Poincar\'e map.  Note that in this case each orientation
of each component of $\{F=\epsilon\}$ give a parameterization
of one of the four circles.

\begin{Cor} \label{c:c1vf}
Let $\gamma$ be simple closed and strictly stable.
There exists $\epsilon>0$ such that if $X_-$ is the vector
field defined on $T_{\epsilon}(\gamma)\setminus \{\gamma\}$ by
$X_-(x)=\tilde \gamma'(0)$ where $\tilde \gamma$ is a noncompact geodesic
with
$\tilde \gamma_-=\gamma$ and $\tilde \gamma (0)=x$, then $X_-$ is $C^1$.
$X_-$ and its (first) derivatives are also continuous functions of the
metric on $M$.
Moreover, there exists a $C^2$
curve $c$ orthogonal to $X_-$ and such that
$\partial c\in \partial T_{\epsilon}(\gamma)\cup \gamma$.
\end{Cor}

The statement that $X_-$ is continuous in the metric makes sense in light of
Lemma \ref{l:vary}.
The vector field $X_-$ is also locally defined and $C^1$
with respect to the metric in a neighborhood of a point $\ell(t_0)$ if
$\ell_- =\gamma$, $\gamma$ is simple closed and
strictly stable, and if the vector field $J_-$
along $\ell$ has $J_-(t_0)\ne 0$.

\begin{proof}
(of Corollary \ref{c:c1vf}).
The continuous dependence of the derivatives of $X_-$ upon the metric
follows from (in order) the continuity of the geodesic flow in the metric;
the fact that the Poincar\'e map ${\mathcal P}$ is $C^1$, with derivatives
depending continuously upon the metric; the fact that the stable manifold
$Y$ is $C^1$, with derivatives depending continuously upon the metric.  It
can also be seen more directly using only the continuity of the geodesic
flow and general dynamic properties of the flow near $\gamma$.
\end{proof}

The next corollary is central to what follows.
A noncompact leaf $\ell$ in a geodesic lamination with $\ell_+$, $\ell_-$
simple closed and strictly stable geodesics
lies in the intersection of
two circles of geodesics, corresponding to its limit
leaves $\ell_{+}$ and $\ell_{-}$.  If $\ell$ is not an
isolated leaf, say $\ell_{i}\in\cL$ with $\ell_{i}(0)
\rightarrow \ell(0)$, then as in Corollary \ref{c:s2graph}
the leaves $\ell_{i}$ have the same limit leaves
and thus also lie in the intersection of the same two
circles of geodesics.  Corollary \ref{c:finitel} can be thought of
as saying that under these assumptions, (which will be satisfied if $\ell$ 
is
a nonisolated leaf in a geodesic lamination on the sphere), the two
circles of geodesics have a common tangent vector at
the point $\ell$; thus the two circles are intersecting
non-transversely.

\begin{Cor} \label{c:finitel}
Let $\cL$ be a geodesic lamination on $M$. If
$\ell,\,\ell_i\in \cL$ are
(distinct)
noncompact, $\ell_+$, $\ell_-$ are strictly stable simple closed
geodesics, and
$\ell_i(0)\to \ell(0)$, then $\ell$ is stable and
there exists a bounded (nontrivial) Jacobi field on $\ell$.
Thus $J_{+}=J_{-}$.
\end{Cor}

\begin{proof}
As in Corollary \ref{c:s2graph} we can assume that
$\ell_{i_{\pm}}=\ell_{\pm}$,
and that $\ell_{i}$ is the normal graph of a function
$u_{i}(t)$ over $\ell$.  Let $\Sigma$ be a local section
near the point $\ell_{+}'(0)$. By Lemma \ref{l:stablemanif}
$\ell$ and $\ell_{i}$
all correspond to points in $\Sigma$ lying on the
stable manifold $Y$.  By Lemma \ref{l:linearna}
$\lim_{i\to\infty} u_{i}/\|(u_{i}(0)),u_{i}'(0))\|$ exists, and is equal to
$ J_{+}$.   By the same reasoning,\\
$\lim_{i\to\infty}  u_{i}/\|(u_{i}(0),u_{i}'(0))\| = J_{-}$.
\end{proof}

\begin{Cor}
For a closed surface with a bumpy metric, the dynamical and 
variational versions of the Kupka-Smale hypothesis are equivalent.  A bumpy 
metric is K--S (by either definition) if and only if for each simple stable
(noncompact) geodesic with closed ends, $J_+ \neq J_-$.
\end{Cor}

\begin{Pro}   \label{p:finitep}
Let $M^2$ be closed surface with a bumpy metric and let $\cL$ be a geodesic
lamination with closed ends.  Then each leaf has finite index 
and $\cL$ has at most
finitely many closed leaves each of which is either isolated or strictly
stable.  Moreover, if $\cL$ has infinitely many leaves, then there exists
a stable noncompact leaf with a (nontrivial) bounded Jacobi field.
\end{Pro}

\begin{proof}
This follows by combining Lemmas \ref{l:nonisolated}, \ref{l:stablengbh}
with Corollaries \ref{c:short} and \ref{c:finitel}.
\end{proof}

There are two different ways of proving our bounds on the Morse indices.
One is to use exclusively the Poincar\'e map and Jacobi fields (this is
the way we will prove Proposition \ref{p:finitep2} below).
The other is to construct positive
supersolutions
of the Jacobi equation.  A particularly simple example
of the second is given in Lemma \ref{l:stablengbh}.   Each 
approach uses the eigenvalue gap, that is that
there are no bounded (nontrivial) Jacobi fields on simple noncompact 
geodesics with closed ends.

\vskip2mm
   Before proving our next result we will need a brief discussion on
convergence of a sequence of simple closed geodesics $\{\gamma_i\}$
in a closed orientable surface $M^2$.  Let $r_0>0$ be sufficiently
small depending only on $\max_M |K|$ and the injectivity radius of $M$.
Fix $x\in M$, then $B_{r_0}(x)\cap \gamma_i$ is the union of disjoint
geodesics segments of length at most $2r_0$ for each $i$.  Note that
any two such that come close to each other are ``almost parallel''.
In fact it follows easily from the equation for geodesics that for each
$i$ there is a coordinate chart $B_{r_0}(x)\to B_{r_0}(0)\subset \RR^2$
such that each component of $B_{r_0}(x)\cap \gamma_i$ is mapped to a line
segment of the form $B_{r_0}(0)\cap (\RR\times \{t\})\subset \RR^2$.  In
this way one can think of each $\gamma_i$ as a geodesic lamination where
the size of the coordinate chart (and the regularity of the maps)
given in the definition of a lamination
is independent of $i$.  Since by the Arzela-Ascoli theorem such a
sequence of coordinate charts is precompact it follows that a
subsequence of the $\gamma_i$'s converge (as a sequence of laminations)
to a geodesic lamination $\cL$.  Implicit in this is that a
sequence of laminations is said to converge if the corresponding coordinate
charts converge and the local transversals converge as closed subsets of
$\RR$ in the Hausdorff
topology.

It follows from this discussion that if $M^2$ is closed with a bumpy metric
and $\{\gamma_i\}$ is a sequence of simple closed geodesics,
then after passing
to a subsequence we may assume that $\gamma_i\to \cL$, where $\cL$ is a
geodesic lamination.  Suppose that $\cL$ has closed ends and 
let $\{\eta_j\}$ be the finitely many closed leaves of
$\cL$ and let $\epsilon>0$ be sufficiently small so that \eqr{e:eprest}
holds.
Let $S$ be given by the remarks surrounding
Corollary \ref{c:short}.   Note that in this case where $\gamma_i\to \cL$,
then the support of $\cL$ is connected and each $\eta_j$ is a nonisolated
leaf, hence strictly stable.  If in addition $M^2$ does not have a
noncompact simple geodesic with a (nontrivial) bounded Jacobi field,
then by Proposition \ref{p:finitep} we can let
$\{\ell_k\}_{k=1,\cdots,m+n}$ be the finitely many noncompact leaves
of $\cL$ ordered so that the first $m$ are the stable leaves.   It follows
that for $i$ sufficiently large each $\gamma_i$ can be decomposed into
pieces that spiral very tightly around one of the $\eta_j$'s and pieces
that are very close and almost parallel to one of the
$(\ell_k\setminus S)$'s.
Note that although there is no a priori bound for how many pieces that
circle one of the $\eta_j$'s or are almost parallel to one of the
$(\ell_k\setminus S)$'s for $k\leq m$ and a given $i$,
only one piece can be very close to an
unstable leaf $\ell_k$.  Had this last claim not been the case
then we would get a contradiction by writing one of the two
(disjoint) pieces as a graph over the other and arguing as in
Corollary \ref{c:finitel} to get a positive Jacobi field.

\begin{Pro}   \label{p:finitep2}
Let $M^2$ be a closed surface with a bumpy Kupka-Smale metric.
If $\gamma_i$ is a collection of simple closed geodesics and every
limit of $\{\gamma_i\}$ is a lamination with closed ends, 
then there is a uniform bound for the Morse indices of $\{\gamma_i\}$.
\end{Pro}

We say that a bumpy metric on
$M^2$ is a bumpy Kupka-Smale metric (or bumpy KS-metric) 
if for each simple stable
(noncompact)
geodesic with closed ends 
in $M$, any bounded (normal) Jacobi field vanishes
identically.
Note that if $M^2$ is a closed surface with a bumpy metric,
then by the results above 
the metric is KS if and only if for each
simple
stable (noncompact) geodesic with closed ends, $J_{-}\ne J_{+}$. 

The proof of Proposition \ref{p:finitep2} will be by contradiction.  We will
assume that $\gamma_i$ is a sequence of simple closed geodesics in
a fixed metric and show that the Morse index of these is uniformly bounded.
A limit of such a sequence is a geodesic lamination.  The relative
Poincar\'{e} map of one of the $\gamma_i$'s (which will tell us how many
zeroes a Jacobi field can have) near the limit lamination will be pieced
together from pieces taken from the leaves of the foliation.  Thus we will
consider first the relative Poincar\'{e} along such a leaf $\gamma$.
The next three lemmas examine the three cases: closed leaves $\eta_j$,
noncompact but stable leaves $\ell_k$ ($k\leq m$) and noncompact unstable
leaves $\ell_k$ ($k>m$).  We need
to watch the image of a fixed vector under the relative Poincar\'{e} map to
see how many times it can cross the $y$-axis.  Recall that such crossings
are always transverse and clockwise.  In order to prevent the corresponding
Jacobi field from having more than one zero we will try to trap it in the
right half plane after it crosses the positive $y$-axis.

First let $\eta$ be simple closed and
strictly stable. Then using the fact that no Jacobi
field has $2$ zeroes, it is not difficult to see that the Jacobi fields
$J_{\pm}$ have no zeroes; thus the vectors
$X_{\pm}(t)=(J_{\pm}(t),J_{\pm}'(t)) \in \RR^2$ never lie along the 
$y$-axis.
If (according to our
convention above) $J_{\pm}(t)>0$
for all $ t$, then $X_-$ and $X_+$ both lie to the right of the $y$-axis, 
and
$X_-$ lies between  the positive $y$-axis and $X_+$.
$P_{t,s}X_{\pm}(s)=X_{\pm}(t)$, and $P_{t,s}$ preserves the four quadrants
defined by $\pm X_{\pm}$.  $P_{t+s_{\eta},t}$ has eigenvectors $X_{-}(t),
X_{+}(t)$ with eigenvalues $1/\lambda,\lambda$ $(0<\lambda<1)$ (the
eigenvalues are independent of $t$) thus the directions of the vectors
$X_{\pm}(t)$ are periodic in $t$. $P_{t+s_{\eta},t}$ fixes these directions
and pushes vectors in the four
``quadrants" away from the $X_+$-direction and toward the $X_-$-direction.
Thus the region between the positive $y$-axis and the vector $X_+$ is a
``trap" from which the future orbit of a vector under $P$ cannot escape.
Let $X_0(t)$ be a (unit) vector in $\RR^2$ halfway between $X_{\pm}(t)$.

\begin{Lem}  \label{l:dkn1}
For each strictly stable simple closed
geodesic $\eta$, there exist $\epsilon$, $H>0$ so that
any
geodesic which is the normal graph
over $\eta|_{[a,b]}$,with $b-a>H$, of a positive function $<\epsilon$,
has the property that if $J$ is the Jacobi field with
$(J(a),J'(a))=X_{0}(a)$, then $J$ has no zero in $[a,b]$, and
$(J(b),J'(b))$ lies above $X_{0}(b)$ (in the right half-plane).
(Here both the geodesic and its Jacobi field have been reparameterized as
graphs.)
\end{Lem}

\begin{proof}
This follows from the fact that by \eqr{e:pigamma} the
curvature $K$ (as a function
of arclength), and thus the relative Poincar\'{e} map
$P_{t+s_{\eta},t}$ for the nearby geodesic will be close to that of
$\eta$, and the fact that the latter has eigenvalues $\lambda$, $1/\lambda$.
This hyperbolicity means that we do not need a bound for $H$.
\end{proof}

Lemma \ref{l:dkn1} also follows from Lemma \ref{l:stablengbh}.

\begin{Cor}  \label{c:dkn2}
If $\cL$ is a geodesic lamination with closed ends on a closed surface
with a bumpy metric, then each
leaf has finite index.
\end{Cor}

If $\cL$ is a geodesic lamination with closed ends on $M^2$ with a 
bumpy metric, then we 
can
find an $\epsilon$ and an $L$ which will work for each closed leaf.  We will
assume this in the next two lemmas.

Now let $\ell$ be be noncompact and stable, with no (nontrivial)
bounded Jacobi field.  Assume $\ell_-$ and $\ell_+$ are strictly stable.
The vectors $X_{\pm}(t)=(J_{\pm}(t),J_{\pm}'(t))$ never cross the
$y$-axis, and as above $X_-$ lies between the positive $y$-axis and $X_+$.
$P_{t,s}$ preserves $\pm X_{\pm}$ and the four quadrants. Let $P_{+}(t)$ be
the (nonrelative; ordinary) Poincar\'{e} map for $\ell_+$, with
eigenvectors $V_{\pm}^{+}(t)$.  For $T$ large there is a map
$\phi:[T,\infty)\to \RR$ so that $\ell_{+}(\phi(t))$ is a
(correctly parameterized) normal graph over $\ell(t)$ on $[T,\infty)$.
For fixed $s$,
\begin{equation}
\lim_{t\rightarrow\infty}\|P_{t+s,t}(\ell)-
P_{\phi(t)+s,\phi(t)}(\ell_+)\|=0\, .
\end{equation}
Thus in particular when $s=s_{\ell_{+}}$, the period of $\ell_{+}$,
$P_{t+s,t}(\ell)$ will have eigenvectors close to $V_{\pm}$ and
eigenvalues close to those of $P_{+}(\phi(t))$.  Given
$\epsilon>0,\mu>0$, for $t$ sufficiently large the image under
$P_{t,s}(\ell)$ of any vector not within an angle $\mu$ of $\pm X_{+}(s)$
(in $\RR^2$ ``at the time $s$'') will lie at an angle $<\epsilon$ of
$V_{-}(\phi(t+s))$.  From this it follows (though different arguments
are needed for the two cases) that
\begin{equation}
   \lim_{t\to \infty}\|\frac{X_{\pm}(t)}{\|X_{\pm}(t)\|}
  -\frac{V_{\pm}(\phi(t))}{\|V_{\pm}(\phi(t))\|}\|=0\, .
\end{equation}
(As $t$ gets large, $V_{-}$ will ``soak up" all vectors except $V_{+}$,
including $X_{-}$. $X_{+}$ does not get soaked up, and must therefore
$=V_{+}$.)  Thus as $\ell$ spirals toward $\ell_+$, the basis
$(X_-,X_+)$ approaches the basis $(V_-,V_+)$, so that the two relative
Poincar\'{e} maps can be glued together.  Similarly, as
$t\rightarrow-\infty$, the basis $(X_-(t),X_+(t))$ approaches a basis
$(U_-(\psi(t)),U_+(\psi(t)))$ of eigenvectors of $P_{-}(\psi(t))$ for an
appropriate reparameterization $\psi$.

\begin{Lem}   \label{l:dkn3}
    For each such (parameterized) geodesic $\ell$,  there is an interval
$[a,b]$ and a $\delta>0$ so that, for each geodesic $\gamma$ which is the
normal graph over $\ell|_{[a,b]}$ of a positive function $<\delta$,\\
  1) $\gamma|_{[a-L,a]}$ and $\gamma|_{[b,b+L]}$ are normal graphs of 
positive
functions $<\epsilon$ over closed leaves as in Lemma \ref{l:dkn1}.\\
2) If $J$ is
the Jacobi field with $(J(a),J'(a))=(U_{0}(\psi(a))$ (a unit
vector midway between $U_+(\psi(a))$ and $U_-(\psi(a))$), then $J$ has no
zero in $[a,b]$, and $(J(b),J'(b))$ lies above $V_{0}(\phi(b))$
(in the right half-plane).\\
  (Here both $\gamma$ and its Jacobi field have
been
reparameterized as graphs.)
\end{Lem}

\begin{proof}
This is similar to the previous lemma.
\end{proof}

\begin{Lem}   \label{l:dkn4}
Let $\ell$ be a noncompact, unstable leaf; assume $\ell_{\pm}$ are strictly
stable simple closed geodesics.  There is an interval $[a,b]$,
a $\delta>0$, and $N\in\ZZ$
  so that, for each
geodesic $\gamma$ which is the normal graph over $\ell |_{[a,b]}$ of a
positive function $<\delta$,\\
  1) $\gamma|_{[a-L,a]}$ and
$\gamma|_{[b,b+L]}$ are normal graphs of positive functions $<\epsilon$
over closed leaves as in Lemma \ref{l:dkn1}.\\
2) If $J$ is a Jacobi field along
$\gamma|_{[a,b]}$, then $J$ has at most $N$ zeroes.
\end{Lem}

We can find $a$, $b$ and $\delta$ which work simultaneously for all the
noncompact leaves.

\begin{proof}
(of Proposition \ref{p:finitep2}).
Suppose now that there were no bound for the Morse index of all simple
closed geodesics on $M$; it follows easily that there exists a sequence
$\{\gamma_i\}$
of simple closed geodesics with index $\to \infty$ and so that $\gamma_i\to
\cL$, where $\cL$ is a geodesic lamination.  Let
$\{\eta_j\}$, $\{\ell_k\}_{k=1,\cdots,m+n}$ be the leaves of $\cL$ as above.
Let $H$, $\epsilon$, $\delta$, $a$, $b$ be as above.
If $i$ is sufficiently large
$\gamma_i$ will consist of a union of pieces, each of which lies within
$\epsilon$ of a closed leaf $\eta_j$ for a time $>H$, or which (after
reparameterization) is a normal graph of magnitude $<\delta$ over some
$\ell_k |{[a,b]}$.  We can assume that only one piece is a normal graph over
$\ell_k$ if $k>m$ (by convention $\ell_k$ is unstable for $k>m$).
Let $J$ be a Jacobi field
along $\gamma=\gamma_i$. We claim that $J$ can have at most $n(N+1)$ zeroes.
By Lemmas \ref{l:dkn1}, \ref{l:dkn3}, if $\gamma_{i}|_{[c,d]}$
is a union of pieces as
above, but only using the closed leaves $\eta_j$ and stable noncompact
leaves, then $J|_{[c,d]}$ can have at most one zero.  Once the vector
$(J(t),J'(t))$ crosses the (say) positive $y$-axis, it will be
trapped in the right half-plane, preventing $J$ from having another zero.
The claim follows using Lemma \ref{l:dkn4}, and thus the proposition.
\end{proof}

\section{Genericity of bumpy Kupka-Smale metrics}  \label{s:s3}

 We begin with some general comments on the Kupka-Smale hypothesis.

 A vector field on manifold is {\em Kupka-Smale} (see \cite{PW}) if
 (1) all closed orbits are hyperbolic, that is, their Poincar\'e maps do not
 have any eigenvalue of modulus $1$, and (2) stable and unstable manifolds
 of closed orbits intersect transversely. Note without hyperbolicity there is
 in general no manifold structure on the invariant sets, so (2) does not make
 sense without (1).  The Kupka-Smale Theorem  states that Kupka-Smale vector
 fields are generic among $C^r$-vector fields, for $r\geq 1$.  Our first
 (dynamical version) definition of KS {\em metric} has in mind the 
 Kupka-Smale condition on the vector field generating the geodesic flow on
 the unit tangent bundle of $M$. There are two differences:  First, we are
 interested only in simple geodesics. Second, we only insist upon
 hyperbolicity for closed geodesics whose Poincar\'{e} maps have real
 eigenvalues.  Note this includes all stable closed geodesics and all closed
 ends of simple geodesics, so the manifolds in (2) are indeed manifolds.
 Since a Poincar\'{e} map coming from the geodesic flow on a surface always
 has determinant $1$, if the eigenvalues are not real there is no hope of
 pushing them off the unit circle by a small change in metric.  Thus in the
 sense of dynamics of simple geodesics our KS-metric condition is the most
 one could ask of a generic metric.  A careful reader can check
 that this condition alone is enough to prove Theorems \ref{t:geolama} and
 \ref{t:geolam}.  

 The variational version of the KS-metric condition is also interesting.
 Condition (1) is a weak version of the bumpy metric condition.  It says
 that, to first order, $\gamma$ does not lie locally in a foliation of $M$ by
 simple closed geodesics. Condition (2) is an analog of the bumpy condition
 for 
{\em noncompact} geodesics, and appears to be independent of the bumpy metric
 condition.   

In this
section we will show
that on $M^2$ bumpy KS-metrics are generic.
Here is the idea:  A
simple noncompact geodesic $\gamma$  is an intersection
point of two circles
of geodesics spiraling toward $\gamma_{-}$ and $\gamma_
{+}$.  The
intersection will be transverse if the two circles have
different tangent
vectors at $\gamma$, i.e., if $J_{-}\neq J_{+}$.  We
will show how to deform
the metric so as to make the two circles intersect
transversely.

Fix a metric $g$ on $M$ and suppose that
$\gamma_1$, $\gamma_2\subset M$ are strictly stable simple closed 
geodesics
(where $\gamma_1=\gamma_2$ is allowed).  Let $\Gamma_{-}
$ be the noncompact
(unit speed) geodesics $\gamma$ with $\gamma_-=\gamma_1
$, and $\Gamma_{+}$
those with $\gamma_+=\gamma_2$.
If $\gamma\in \Gamma_-$, then we let
$F_-:(-\epsilon,\epsilon)\times \RR\to M$ be a
(nontrivial)
geodesic variation of $\gamma$ so that $F_-(s,\cdot)\in
\Gamma_-$,
$F_-(0,\cdot)=\gamma$, and
$g(\frac{\partial F_-}{\partial s},\frac{\partial F_-}
{\partial t})|_{s=0}=0$.
Likewise if $\gamma\in \Gamma_+$, then we let $F_+$ be a
(nontrivial) geodesic variation of $\gamma$ consisting
of geodesics
asymptotic to $\gamma_2$.
(Note that by Section \ref{s:s2} these variations are
essentially unique; we can also assume that $\frac
{\partial}{\partial
s}F_{\pm}(0,\cdot)= J_{\pm}(\gamma)$ and $\frac
{\partial^2}{\partial
s\partial t}F_{\pm}(0,\cdot)=J_{\pm}'(\gamma)$). 
We say that $\Gamma_-$ and $\Gamma_+$
intersect transversally at $\gamma$ if the two curves
representing
$\Gamma_{-}$ and $\Gamma_{+}$ in a local section
$\Sigma$ at $\gamma$ are transverse,
i.e., if given a curve
$\nu:(-\epsilon,\epsilon)\to M$ with $\nu (0)=\gamma
(t_0)$ for some $t_0$
and $\nu '(0)$ transverse to $\gamma '(0)$, the curves
$\frac{\partial F_-}{\partial t} |\nu$ ($=X_{-}|\nu$)
and
$\frac{\partial F_+}{\partial t} | \nu$ ($=X_{+}|\nu$)
in $\Pi ^{-1}(\nu)\subset T_1M$
are (well defined and) transverse at $0$. (Note that
$\frac{\partial F_-}{\partial t}(0,0)=\frac{\partial
F_+}{\partial t}
(0,0)=\gamma '(t_0 )$.)
Now a tangent vector to the curve $\frac{\partial F_-}
{\partial t}|\nu$ at
$s=0$ is given by
$(\nu_{\ast}(\frac{d}{d\tau}) F_{-}(0,t_{0}),\nu_{\ast}
(\frac{d}{d\tau})
\frac{\partial F_{-}}{\partial t}(0,t_{0}))$, which is proportional to
$(J_{-}(t_{0})\,\nn+\alpha \frac {\partial}{\partial t},J_{-}'(t_{0})\,\nn)$, 
where $\alpha$ is the slope
of the tangent vector
to $\nu$ in the $(s,t)$ coordinates, and similarly for
$F_{+},J_{+}.$
Thus transversality means that the vectors $(J_{-}(t_
{0}),J_{-}'(t_{0}))$
and $(J_{+}(t_{0}),J_{+}'(t_{0}))$ are not parallel.
Since the Jacobi
fields
are determined uniquely by their and their derivatives values at $t=t_{0}$,
by Section \ref{s:s2}
transversality
of $\Gamma_-$ and $\Gamma_+$ at $\gamma$ is equivalent
to the fact that $\gamma$
has no bounded nontrivial Jacobi field.
To prove that the set of bumpy KS-metrics on
$M^2$
is residual it suffices
therefore (by Lemma \ref{l:vary})
to show that a residual set of metrics on $M^2$
consists of bumpy metrics with the property that, for
each pair
of strictly stable simple closed
geodesics $\gamma_1$ and $\gamma_2$, $\Gamma_-$ and
$\Gamma_+$ intersect
transversally.

The rest of this section is devoted to the proof of:

\begin{Thm}  \label{t:strictlybumpy}
On a closed surface the set of bumpy KS-metrics contain a
residual set.
\end{Thm}

\begin{proof}
(of Theorems \ref{t:geolama} and \ref{t:geolam}).
Theorems \ref{t:geolama}, \ref{t:geolam} follows by
combining
Theorem \ref{t:strictlybumpy} with Propositions \ref
{p:finitep},
\ref{p:finitep2}, respectively.
\end{proof}

In Lemmas \ref{l:cuta}, \ref{l:cutandpaste}
below we let $M^2$ be a closed
surface with metric $g$ and $c:[-\epsilon,r_0+\epsilon]
\to M$
be a simple $C^{m+2}$ curve
parameterized by arclength.
Let $\nn_c$ be the unit normal (so $\nn_c\in C^{m+1}$)
and $\Phi:[-\epsilon,s_0+\epsilon]\times [0,\epsilon]
\to M$ given
by $\Phi(s,t)=\exp_{c(s)}(t\,\nn_{c(s)})$ so $\Phi^{-1}
$ are
geodesic normal coordinates
in a neighborhood of $c$.
In these the metric is of the form $f^2(s,t)ds^2+dt^2$
where $f\in C^m$.
(Note that in any metric of this form, the curves
$s=$constant are minimal
and thus geodesic.)

We will give a simple direct argument to show that the
set of strictly
bumpy metrics on $M$ are dense.  The following
deformation lemma
will be needed to show that on a surface
for a dense set of metrics certain geodesic variations
intersect transversally.  The lemma allows us to alter
the geodesic flow in
a controlled way by altering the metric.

\begin{Lem}   \label{l:cuta}
Let $M^2$, $c$, $\Phi$ be as above. Let $c_w:[-
\epsilon,s_0+\epsilon]\to M$
be the curve
$c_w(s)=\Phi (s,\epsilon-w\,s)$.
There exists a $1$-parameter family of $C^m$ metrics
$g_w$ ($w\in (-\delta,\delta)$)
such that $g_0=g$, each $g_w=g$ on
$M\setminus \Phi([-\epsilon,s_0-\epsilon]\times [0,3
\epsilon/4])$,
$\Phi_w (s,\epsilon)=c_{w}(s)$ for $s\in [0,s_0]$, and
$\frac{\partial}{\partial t}\Phi_{w}(s,\epsilon)$
points perpendicular to
the curve $c_{w}(s)$.
(Here $\Phi_w^{-1}$
are geodesic normal coordinates in a neighborhood of
$c$ in the metric
$g_w$.) That is, in the $g_w$ metric, the family of
geodesics which enter the
box $\Phi$ perpendicular to the curve $c(s)$ exit
perpendicular to the curve
$c_{w}(s)$.
\end{Lem}

\begin{proof}
Fix $w>0$ sufficiently small.  Let $\Psi_w^{-1}$
be geodesic normal coordinates in a neighborhood of
$c_w$
parameterized so that
$\Psi_w(s,t):[-\epsilon,s_0+\epsilon]\times
[0,\epsilon]\to M$,
and $\Psi_w(s,\cdot)$ are
geodesics moving away from $c$ and ending up on $c_w$
orthogonal to
$c_w$.
In particular $\Psi_w(s,\epsilon)=c_w(s)$.  In these
coordinates the metric can be written as $f_w^2(s,t)
ds^2+dt^2$. Let
$\eta:[0,\epsilon)\to [0,1]$ be a smooth cutoff
function with
$\eta |[0,\epsilon/2]=1$.  Then
$\Phi^{-1}_{\eta,w}(s,t)=(\eta(t)\,\Phi^{-1}(s,t)
+(1-\eta (t))\,\Psi_w^{-1})(s,t)$ is a
diffeomorphism and give therefore local coordinates
$(s,t)$.
In these define a
metric by
$\hat{g}_w(s,t)=(\eta (t) f(s,t)+(1-\eta (t))f_w(s,t))
^2ds^2+dt^2$.
Finally, let
$\phi\in C_0^{\infty}(-\epsilon,s_0+\epsilon)$ with $0
\leq \phi\leq 1$,
$\phi | [0,s_0]=1$, and set
$g_w (s,t)=\phi (s)\,\hat g_w (s,t)+(1-\phi (s))\,g
(s,t)$.  It
is easy to see that this gives a
$1$-parameter family with the desired properties.
\end{proof}

We will use this deformation to make the family $\Gamma_
{-}$ (locally)
transverse to $\Gamma_{+}$.  The metric will be
deformed in a rectangle
$\Phi$ to change the  family $\Gamma_{-}$ as it moves
through the rectangle,
(roughly speaking) before it meets the family $\Gamma_
{+}$ at the top of the
rectangle.
To see the effect of this deformation of the metric on
the image of the
family of geodesics $\Gamma_{-}$ in the local section
$\Sigma$ given by the
(fixed) curve $c_{0}$ at the top of the rectangle, we
will need to know the
angle at which
these geodesics cross the
curve $c_{0}$.  The image of the family in $\Sigma$ is
(in appropriate
coordinates) the graph of the
crossing angle as a function of arclength along $c_0$.
Lemma \ref{l:cutandpaste} below
begins by showing that the deformation of the metric
given in
Lemma \ref{l:cuta}
moves this graph (and thus the curve which is the image
of the family
$\Gamma_-$ in $\Sigma$) off itself.
The following version of Sard's
theorem says that if we can move a curve in the plane
off itself, we can
make it transverse to a second curve.  (Here we say that
the intersection of the image of $h$ and the graph of
$f$ is transverse if, for every $(x,y)$ and $t$ with
$h(t)=(x,y)$, $h'(t)=(\frac{dx}{dt},\frac{dy}{dt})$ is
transverse to
$(1,f'(x))$.)

\begin{Lem}  \label{l:sard}
Let $f_w(s)$ be $C^1$ functions, where
$s\in [0,1]$ and $w\in (-\delta,\delta)$.  Let
$h:\RR\rightarrow\RR^2$ be
$C^1$.  If
$|\frac{\partial f_w}{\partial w}|_{w=0}$ is
nonvanishing,
then there exists a sequence $w_i\to 0$ such
that the curves $(s,f_{w_i}(s))$ are
transverse to $h$.
\end{Lem}

\begin{proof}
By a $C^1$ change of coordinates, we can assume the
functions $f_{w}(s)$ are
constant.  In these coordinates transversality of $h$
and $f_{w}$ is
equivalent to $f_w$ being a regular value of $\pi_2
\circ h$, where $\pi_2$
is projection onto the second factor in $\RR^2$.
By Sard's theorem the claim easily follows.
\end{proof}

\begin{Lem}  \label{l:cutandpaste}
Let $M^2$, $c$, $\Phi$ be as above with $m=0$.
If $h:[0,1]\to T_1M|\Phi (\cdot,\epsilon)$ is a
$C^1$ curve, then there exists a sequence of metrics
$g_i\to g$ with $g_i=g$
on $M\setminus \Phi([-\epsilon,s_0+\epsilon]\times [0,3
\epsilon/4])$
and such that in any $g_i$,  $\frac{\partial \Phi_i}
{\partial t}$ intersects
$h$ transversally along $\Phi (\cdot,\epsilon)|[0,s_0]
$.  Here $\Phi_i^{-1}$
are geodesic normal coordinates in a neighborhood of $c$
in the metric $g_i$.
\end{Lem}

\begin{proof}
Assume first that $c$ is actually $C^{\infty}$.  As
above let
$\Phi^{-1}=(s,t)$
be geodesic normal coordinates in a neighborhood of $c$
so that we can
in particular think of $t$ as a function on this
neighborhood.
Moreover, in this
neighborhood the metric is $f^2(s,t)ds^2+dt^2$.  Let
$g_w$,
$\Phi_w$ be
given by Lemma \ref{l:cuta} and set
$h_w=\frac{\partial \Phi_s}{\partial t}|\Phi
(\cdot,\epsilon)$.  It
follows from Lemma \ref{l:cuta} that
$\left| \frac{\partial }{\partial w}_{|w=0}h_{w}\right|$
is nonvanishing for $s\in [0,s_0]$.
Namely, it is easy to see that $h_{w}$ is $C^1$ so
we need only check that the derivative is
nonvanishing.  To see this let
$\gamma_w:[\theta_0,\theta_w]\to M$ be a (unit speed)
geodesic (in the metric $g$) with
$\gamma_w (\theta_0)\in\{\Phi (s,\epsilon-ws)\}$,
$\gamma_w'(\theta_0)$
orthogonal to  $\{\Phi (s,\epsilon-ws)\}$, and
$\gamma_w (\theta_w)\in \{\Phi (s,\epsilon)\}$.  Set
$u(\theta)= (t\circ \gamma_w)'(\theta)=g(
\gamma_w'(\theta), \nabla t)$,
then
$u'(\theta)=\Hess_t(\gamma_w',\gamma_w')=(1-u^2(\theta))
\,\frac{f'}{f}$.
Hence, $|(\log [(1+u)/(1-u)])'|\leq C$, where $C=C(g)$
is a constant.
In particular,
\begin{equation}  \label{e:cutandpaste1}
\left| \frac{(1+u(\theta_w))}{(1-u(\theta_w))}
                  \frac{(1-u(\theta_0))}{(1+u(\theta_0))}
\right|
                              \leq \exp (C\,|\theta_w-
\theta_0|)
\, .
\end{equation}
Moreover, it is easy to see that for some
$\alpha=\alpha (g)$
and some $\beta=\beta (g)>0$
\begin{equation}  \label{e:cutandpaste2}
|\theta_w-\theta_0|\leq \alpha\, w\text{ and }|u
(\theta_0)|\geq \beta\, w\,
.
\end{equation}
Combining \eqr{e:cutandpaste1} with \eqr
{e:cutandpaste2} we
conclude that for some $\beta'=\beta'(g)>0$
\begin{equation}  \label{e:cutandpaste3}
|u(\theta_w)|\geq \beta'\, w\, .
\end{equation}
Since $g_w=g$ on
$M\setminus \Phi([-\epsilon,s_0+\epsilon]\times [0,3
\epsilon/4])$
it follows easily from \eqr{e:cutandpaste3} that
$\left| \frac{\partial }{\partial w}_{|w=0}h_w\right|$
is nonvanishing for $s\in [0,s_0]$
and the lemma now follows from Lemma \ref{l:sard}.

In the general case where $c$ is only $C^2$
let $c_j$ be a sequence of $C^{\infty}$ curves in $M$
with
$|c-c_j|_{C^2}\to 0$.  It follows easily from the
continuity of the
quantities
involved that for $j$ sufficiently large (but fixed) if
$g_{w,j}$ and
$\Phi_{w,j}$
are given by Lemma \ref{l:cuta} with respect to $c_j$
and
$h_{w,j}=\frac{\partial \Phi_{w,j}}{\partial r}|\Phi
(\cdot,\epsilon)$, then
$\left| \frac{\partial }{\partial w}_{|w=0}h_{w,j}
\right|$
is nonvanishing.
The lemma now easily follows from Lemma \ref{l:sard}.
\end{proof}

Let again $M^2$ be a closed surface with a metric $g$
and suppose that
$\gamma_1$, $\gamma_2$ are strictly
stable simple closed geodesics
(where $\gamma_1=\gamma_2$ is allowed).  Let $\Gamma_{-
,+}^{k}$ be the
$\gamma\in \Gamma_{-}\cap \Gamma_+$ with
$\gamma\setminus T_{\delta}(\gamma_1\cup\gamma_2)$ of
length $\leq k$.

The next result will follow by applying Lemma \ref
{l:cutandpaste} a
finite number of times:

\begin{Lem}  \label{l:denseglk}
Let $M^2$, $g$, $\gamma_1$, and $\gamma_2$ be as
above.  Given $k>0$,
there exists $g_n\to g$ with $g_n=g$ on in a
neighborhood of
$\gamma_1\cup \gamma_2$
and such that for each $g_n$, then $\Gamma_-(g_n)$
and $\Gamma_+(g_n)$ intersect transversally at $\Gamma_
{-,+}^k(g_n)$.
\end{Lem}

\begin{proof}
Fix $\delta$ small but positive.
Let $\sigma_{1}:\SS^1\rightarrow M$ be a simple closed
curve in
$T_{\delta/2}(\gamma_{1})$ meeting each geodesic in
$\Gamma_{-}$ exactly
once (Corollary \ref{c:stablemfldf}).  Parameterize
$\Gamma_{-}$ as a map
$F_{-}:\SS^{1}\times \RR \to T_{1}M$ with $F_{-}(s,0)
=\sigma_{1}(s)$.
Similarly let $\sigma_{2} :\SS^1\to T_{\delta /2}
(\gamma_{2})$ meet
each geodesic in $\Gamma_{+}$ exactly once, and let
\begin{equation}
F_{+}:\SS^1 \times \RR \rightarrow T_{1}M \mbox{ with }
F_{+}(s,0)=\sigma_{2}(s)\, .
\end{equation}
(Of course we want $F_{-}\to \gamma_{1}$ as $t\to -
\infty$
and $F_{+}\to\gamma_{2}$ as $t\to \infty.$)
Given $\epsilon$ sufficiently small, and $x\in T_
{\delta/2}(\gamma_{1})$, we
can find (Corollary \ref{c:c1vf}) a unit speed curve
$c:[-2\epsilon,2\epsilon]\to T_{\delta}(\gamma_{1})$
with $c(0)=x$,
and which is everywhere perpendicular to the vector
field $X_{-}$.  Let
$\Phi:[-2\epsilon,2\epsilon]\times [-\epsilon,0]\to M$
be geodesic
normal coordinates, with $\Phi(s,0)=c(s)$ and $\Phi
(s,\cdot)\in\Gamma_{-}$
(when extended) for each $s$.  Given $k$, and
$\gamma\in\Gamma_{-}$, there exist $\epsilon,\mu >0$
(depending only upon
$k$, not on $\gamma$) and $t_{0}\in [-1,0]$ (depending
upon $k$ and
$\gamma$)
so that this box $\Phi$ for the point $x=\gamma(t_{0})$
has the following
properties:\\
a) Geodesics entering the bottom ($t=-\epsilon$) of the
box vertically
(tangent to $s=$ constant)\\
    have never been in the box before.\\
b) Geodesics leaving the top ($t=0$) of the box at an
angle less than $\mu$
of the
vertical travel at
least a distance $5k$ before returning within $\epsilon$
of the box.\\
Let $\hat{c}:[-2\epsilon,2\epsilon]\rightarrow\Pi^{-1}
(c)$ be the lift of
the vector field $X_{-}$: $\hat{c}=X_{-}|_{c}.$  Note
that by a)
above,
$\hat{c}$ is the intersection of the local section
$\Sigma
=\Pi^{-1}(c)\subset T_{1}M$ given by the
transversal $c$, with the image of $F_{-}|_{(-
\infty,0]}.$  Let $D$ be the
disk in $\Sigma$ consisting of all vectors in
$\Pi^{-1}(c)$ making an angle less than $\mu$ with $X_{-
}$.  Consider the
image in $T_{1}M$ of $F_{+}|_{[-k,\infty )}$.  This
``tube'' will intersect
$D$ transversely if $\mu<\pi/2$.  The tube will thus
intersect $D$ in a
union $C_k$ of $C^1$ curves whose total length is
finite.
An arbitrarily small alteration $g_{w}$ of the metric
in the box will not
change the intersection $C_k$ of $D$ with $F_{+}|_{[-
k,\infty )}$ by b)
above,
but will, by Lemma \ref{l:cutandpaste}, make
$\hat{c}|_{[-\epsilon,\epsilon]}$
transverse
to $C_k$.

Given $k$ we can find a finite number of such
boxes $\Phi$ with each
geodesic in $\Gamma_{-}$ intersecting some $c|_{[-
\epsilon,\epsilon]}$. Thus
we will be done if we can make each $\hat{c}|_{-
\epsilon,\epsilon]}$
transverse to $C_k$.  Now
when we alter the metric in the second box it may
change the image in the
section $D$ at the top of the first box, of the forward
tube
$F_{+}|_{[-k,\infty )}$.  However a sufficiently small
change of metric in
the second box will not destroy transversality for the
first box due to the
fact that, as a consequence of Corollary \ref
{c:stablemfldf},
transversality is an open
property.
\end{proof}

Since transversality of the sets $\Gamma_-(g)$
and $\Gamma_+(g)$ along
$\Gamma_{-,+}^k(g)$ for any fixed $k>0$ is an open
property in $g$
(this follows easily from Corollary \ref{c:c1vf}) we
get:

\begin{Cor}  \label{c:openglk}
Let $M^2$, $g$, $\gamma_1$, and $\gamma_2$ be as
above.  If $\Gamma_-(g)$
and $\Gamma_+(g)$ intersect transversally at $\Gamma_{-
,+}^k(g)$, then
there exists an open neighborhood $U$ of $g$ such that
for each $\bar{g}\in U$, $\Gamma_-(\bar{g})$ and
$\Gamma_+(\bar{g})$ intersect
transversally at $\Gamma_{-,+}^{k}(\bar{g})$.
\end{Cor}

\begin{proof}
Suppose not; it follows easily that there exists
$g_i\to g$, $\gamma_{1,i}\to \gamma_1$,
$\gamma_{2,i}\to \gamma_2$, and $\gamma_i\in \Gamma_{-
,+}^k(g_i)$ with
$(\frac{\partial F_{-,i}}{\partial s},
\frac{\partial^2 F_{-,i}}{\partial s\partial t})
=(\frac{\partial F_{+,i}}{\partial s},
\frac{\partial^2 F_{+,i}}{\partial s\partial t})$ at
$\gamma$.  Clearly
$\gamma_i\to \gamma\in \Gamma_{-,+}^k(g)$ and by
Corollary \ref{c:c1vf} it follows easily that
$(\frac{\partial F_{-}}{\partial s},
\frac{\partial^2 F_{-}}{\partial s\partial t})
=(\frac{\partial F_{+}}{\partial s},
\frac{\partial^2 F_{+}}{\partial s\partial t})$ at
$\gamma$ which
is the desired contradiction.
\end{proof}

\begin{proof}
(of Theorem \ref{t:strictlybumpy}).
Fix integers $L$, $k>0$ and let
${\mathcal{G}}_{L,k}$ be the set of metrics $g$ on $M^2$
with the two properties:\\
a) All closed geodesics with length $< L$ are
nondegenerate critical
points.\\
b) If $\gamma_1$, $\gamma_2$ are simple closed strictly
stable geodesics
with length $<L$, then $\Gamma_-(g)$ and $\Gamma_+(g)$
intersect
transversally along $\Gamma_{-,+}^k(g)$.

Combining the fact that the set of bumpy metrics is
dense (see \cite{Ab}, \cite{An})
with Lemmas \ref{l:vary}, \ref{l:denseglk},
and Corollary \ref{c:openglk} we
get that ${\mathcal{G}}_{L,k}$ is open and dense, hence
$\cap_{L,k>0}{\mathcal{G}}_{L,k}$ is residual.
\end{proof}

\section{Geodesic laminations without closed ends; Theorem \ref{t:notclosed}}
\label{s:s4}
  
A {\em train-track} is a one-complex $T$ embedded 
in a surface satisfying conditions of (1) smoothness, 
(2)  nondegeneracy, and 
 (3) geometry.  The definition is quite involved, and probably
  familiar to 
 many readers.  Rather than attempt an abridged, incorrect version, we
  refer 
 the reader to \cite{HaPe} (see p. 4 there) for this and other 
 definitions. 

Let $F$ be a disk with four holes removed, see fig. 3a, and let $N$ be
the topological double of $F$ so $N$ is a closed orientable surface of genus
$4$.  Equip $N$ with a metric with negative curvature so that the
boundary of $F\subset N$ consists of geodesics.  

\medskip
\centerline{\epsfig{figure=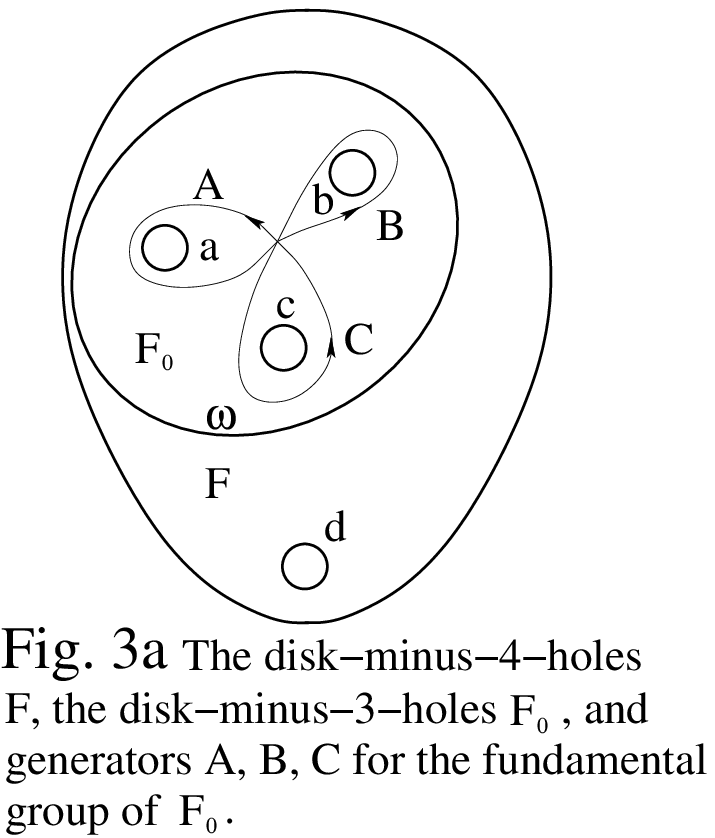}\hspace{1.5cm}
\epsfig{figure=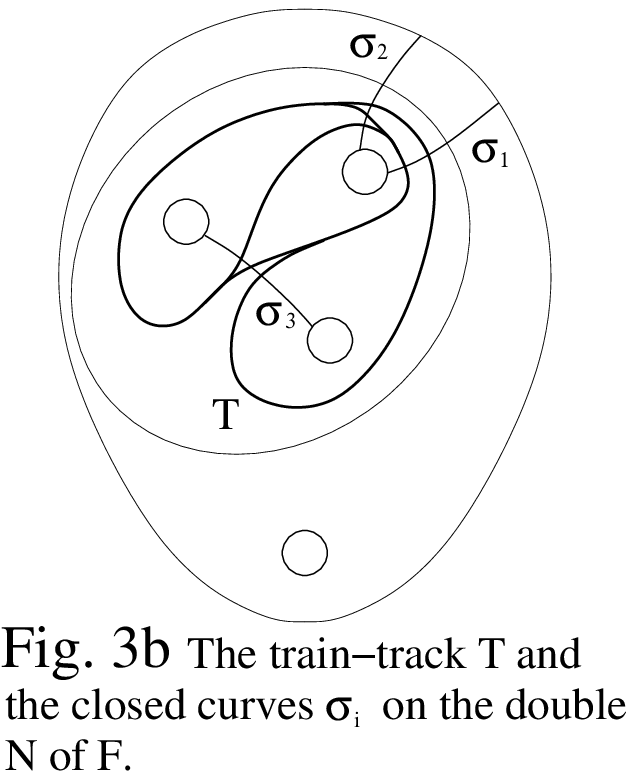}}
\medskip

On $F$ pick 
three of the holes; there is a unique closed geodesic 
$\omega$ on $F$ enclosing these three holes but not the 
fourth.  Let $F_0\subset F$ be the disk bounded by $\omega$, 
with the three holes removed. The fundamental group of $F_0$ is a 
free group on three generators $A$, $B$, $C$, one for each
hole. Fig. 3b  
shows a transversely recurrent train-track $T$ on
$N$. (That $T$ satisfies the ``geometry'' condition for a
train-track is clear from the characterization on the top of p. 6 in
\cite{HaPe}:  The
 complement of $T$ in $N$ is connected, and clearly has no component that is
 an embedded nullgon, monogon, bigon, annulus, or punctured annulus.
 Transverse recurrence can be seen by looking at the closed curves on $N$
 that are the doubles $\sigma_i$ of the curves in fig. 3b.  The reader will
 easily see how to find curves as the $\sigma_i$'s so that each branch of $T$
 meets at least one of the $\sigma_i$, and for each $i$, the complement in
 $N$ of $T\cup \sigma_i$ is either connected, or has two components, one of
 which is an embedded trigon.  Since neither component is an embedded bigon,
 each $\sigma_i$ hits $T$ efficiently (\cite {HaPe}, p. 19). )

\medskip
\centerline{\epsfig{figure=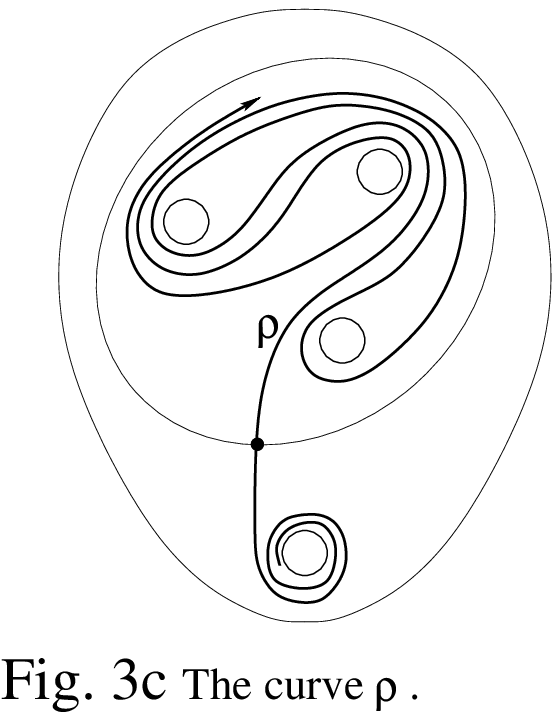}}
\medskip

In fig. 3c we see a simple curve $\rho$ on $F$  that for $t>0$ is
 carried on this train-track.  The bounded homotopy class (that is, we
 only consider homotopies that move points a bounded distance in the
 universal cover) of the forward end of the
 curve $\rho$ determines a (semi-)infinite word 
$BA^{-1}B^{-1}CBAB^{-1}A\cdots$
 in the generators $A$, $B$, $C$.  
Theorem \ref{t:notclosed} will follow from two lemmas:

\begin{Lem} \label{l:A} 
An uncountable number of different words in 
$A$, $B$, $C$ come from different completions of the forward end of $\rho$
as a simple curve carried on the train-track $T$.   
\end{Lem}

\begin{Lem} \label{l:B} 
Each such completion of $\rho$ is bounded-homotopic in
$F$ to a simple geodesic.    
\end{Lem}

\begin{proof} 
(of Lemma \ref{l:A}). 
In fig. 3c each of the holes 
$a$, $b$, $c$ lies in a snake-like cavity bounded 
by a segment of $\rho$.  Two of
 these cavities open out where $\rho$ meets $\omega$; the third ends at the
 forward end of $\rho$.  The retraction of $\rho$ onto $T$ retracts
 the head of the snake to the
 boundary of a monogon (the monogon actually lies in the disk, before the
 holes have been removed to form $F$) whose interior contains 
no other branches of
 $T$, and identifies the two sides of each snake's body, so a simple curve
 carried by $T$ that enters one of these cavities has exactly two ways to
 continue up to homotopy: Clockwise or counterclockwise. These properties
 will persist as we complete the curve $\rho$.  

\medskip
\centerline{\epsfig{figure=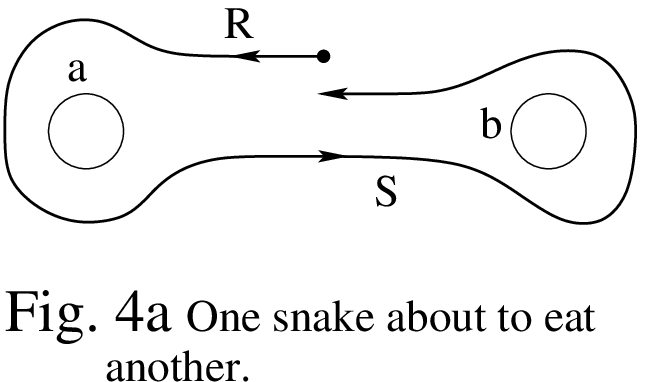}\hspace{1.5cm}
\epsfig{figure=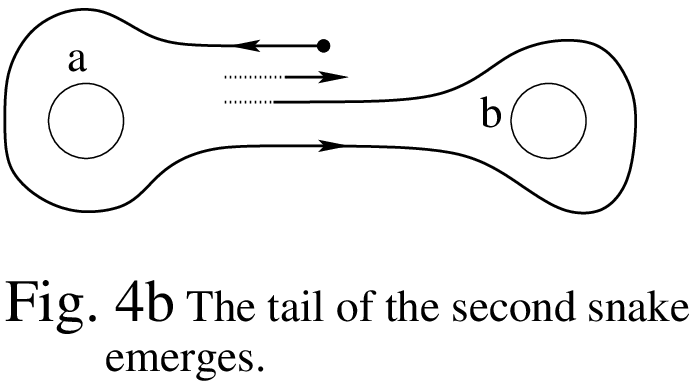}}
\medskip

 At present the forward end of $\rho$ lies at the top of the picture.  It
 will proceed around to the bottom, and enter the cavity containing hole
 $a$.  The reader can easily verify that when one snake eats another as in
 fig. 4a the possibilities are as follows:  Suppose the left snake has
 the word $R$ in the letters $A$, $B$, $C$, and the right 
snake the word $S$.  ($R$
 is conjugate to $A$, and $S$ to $B$.)
\begin{enumerate}
 \item The word $(RS)^kRSR^{-1}(RS)^{-k}$, $k\geq0$.
 \item The word $(RS)^kR(RS)^{-k}$, $k\geq1$.
 \item The infinite word $RSRSRSRSR\cdots$
 \end{enumerate}
 In the first two cases, the forward end emerges as in fig. 4b.  It will
 then proceed to enter the third cavity, and the above possibilities repeat.
 The important things for us are:  There are no ``wrong turns", that is the
 curve can always be completed as a simple curve carried by $T$; and the
 curve never runs out of possibilities, that is there is 
always another choice of
 {\em words} in the future.  It follows that an uncountable number of
 different words results.
 \end{proof}

\begin{proof} 
(of Lemma \ref{l:B}). 
Let $\tilde N$ be the universal cover of $N$ 
and $\tilde N(\infty)$ the sphere (circle) at
 infinity of $\tilde N$; see for instance \cite{Eb}.  
Let $\tilde{\rho}$ be a lift of $\rho$ 
to $\tilde N$, with $p$ the
 point lying over the intersection of $\rho$ with $\omega$.  Let
 $\tilde{\omega}$ be the lift of $\omega$ through $p$, and $\tilde{F}$ be the
 lift of $F$ through $p$.  By proposition 1.5.1 of \cite{HaPe} (that
 proposition is stated in case of constant negative
curvature but the argument extends) 
$\tilde{\rho}$ has two (unique, and distinct) limit points $q,r$ on
 $\tilde N (\infty)$.  Let $\tilde \ell$ be the
 geodesic on $\tilde N$ with these two limit points.  If $\tilde{\tau}$ is a
 geodesic on $\tilde N$ with neither $q$ nor $r$ as a limit point, then the
 topological intersection numbers $\#(\tilde{\rho},\tilde{\tau})$ and
 $\#(\tilde{\ell},\tilde{\tau})$ are well defined and equal.  Because
 $\tilde N$ has
 negative curvature, the number of points of intersection of two geodesics in
 $\tilde N$ is $0$ or $1$.  It follows that $\ell$ lies on 
$F$, and on $F_0$ after
 crossing $\omega$.  Moreover, if $\tilde \sigma$ 
is a geodesic segment, with boundary on $\partial\tilde{F}$, then
 $\#(\tilde{\rho},\tilde{\sigma})=\#(\tilde{\ell},\tilde{\sigma})$.  
The
 reader can easily verify that for each branch of $T$ there is a curve $\tau$
 on $F_0$,  with boundary on $\partial F_0$, that intersects $T$ exactly
 once.  By making $\rho$ stick close to $T$, we can assume that $\tau$
 and $\rho$ intersect once each time (and only when) $\rho$ follows
 the given branch of $T$ (i.e., when the retraction of $\rho$ onto $T$
 includes the branch).  
 Each lift of $\tau$ to $\tilde{F_0}$ that meets $\tilde{\rho}$, say
 at a point $p_i$, is homotopic to a geodesic $\tilde{\sigma}_i$ with boundary
 on $\partial\tilde{F_0}$, and  $\tilde{\sigma}_i$ will 
 intersect $\tilde \ell$ exactly
 once, say at $q_i$. We can now define a homotopy from $\tilde{\rho}$ to
 $\tilde{\ell}$ that starts by taking each point $p_i$ to $q_i$; then ``pull
 tight" to straighten out the curve inbetween the points $q_i$.  Thus
 $\tilde{\rho}$ is homotopic to $\tilde{\ell}$ inside $\tilde{F}$, and
 $\rho$ is homotopic to $\ell$ inside $F$, by a
 homotopy that does not move any point very far.  

 It remains only to see that $\ell$ is simple. Since $\rho$ is simple,
 $\tilde{\rho}$ does not intersect any other lift of $\rho$.  The same
 intersection number argument as before shows that $\tilde{\ell}$ does not
 intersect any other lift of $\ell$, which implies that $\ell$ is
 simple. 
\end{proof}

\begin{proof} (of Theorem \ref{t:notclosed}).  
First we will produce an open set of metrics on $M$ having geodesic
laminations without closed ends.  
Fix a metric on $N$ with constant curvature $-1$ ($F$ is one-half of
$N$) and let $U\subset N$ be an
open neighborhood of $F$.  On the given surface $M$ we can complete
the metric on $U$ to a metric on $M$.  Any nearby metric on $M$ 
will contain a (unique) metric
 surface $F\subset U$ of negative curvature and geodesic boundary that also
 has a completion to a  metric of negative curvature on $N$.

The statement now follows from the lemmas as 
follows:  If $\ell$ is the geodesic in $F$ isotopic to $\rho$, and $\ell$ has 
closed ends,  by Lemma \ref{l:+closed} the word determined by $\rho$
will eventually repeat. But uncountable many of the words determined
by simple completions $\rho$ will be nonrepeating. For such $\rho$ the
closure of $\ell$ will be a lamination of $F$ with nonclosed ends.

Next we find a (smaller) open set of metrics, each  having a geodesic
 lamination without closed ends that is the limit of a sequence of simple
 closed geodesics. 

 Let $F_1$ be a topological disk-minus-6-holes.  As before we assume a metric
 of negative curvature on $F_1$ with geodesic boundary, that extends to a
 metric of negative curvature on the topological double $N_1$ of $F_1$.  Let
 $\omega$ be a geodesic enclosing holes $a$, $b$, $c$ but not holes
 $d$, $e$, $f$.  Put
 the train-track $T$ inside $\omega$ as before, and a similar train-track $S$
 around the holes $e$, $f$, $g$.   
As before we can find a simple curve $\rho$ on
 $F_1$ that crosses $\omega$ once, and whose bounded homotopy class
 determines  a (doubly) infinite word starting in $d$, $e$, $f$ and ending in
 $a$, $b$, $c$, and  not repeating at either end.  We will 
next describe a sequence
 of simple closed curves $\rho_k$ on $F_1$,  which ``approximate" $\rho$. The
 curve $\rho_k$ starts at a point on $\omega$, and follows the path of $\rho$
 in and out of $k-1$ cavities.  
 In the $k$'th cavity, the curve $\rho_k$ will turn
 around counterclockwise, and then stick close to the (just laid out) strand
 of $\rho_k$  on its left, until it retraces its steps back to the
 intersection with $\omega$.  The other half of the path of $\rho_k$ does the
 same thing on the other side, then closes up.  There will be a
 simple closed geodesic $\gamma_k$ in the free homotopy class of
 $\rho_k$.  Take lifts $\tilde{\rho}_k$ through
 a fixed point $p$ lying above $\omega$, and corresponding lifts 
$\tilde{\gamma}_k$. The
 curve $\rho$ will cross (in order) a sequence $\tilde{\sigma_j}$
 of lifts of (a finite number of) geodesic segments $\sigma_j$ on $F_1$, with
 boundary on $\partial F_1$, that keep track of which choice was made at each
 branch.  The same
 will be true of the geodesics $\gamma_k$, as long as $\rho_k$ follows
 $\rho$.  Now fix $J$.  By the Arzela--Ascoli theorem, 
for each $J$ any limit lamination of the sequence $\gamma_k$ will
 contain a geodesic $\ell_J$ with a lift that has the same intersection with
 the $\sigma_j, \|j\|\leq J$, as $\rho$.  Finally, the closure of the
 set of geodesics $\ell_J$ in the unit tangent bundle 
contains a geodesic $\ell$ with a
 lift that intersects all of (and only) the segments $\tilde{\sigma_j}$
 intersecting $\rho$, in order.  This geodesic does not have closed ends,
 since its ``word" is the same as that of $\rho$ (note their lifts are
 bounded--homotopic), and hence
 nonrepeating.  By construction any limit lamination of the 
sequence $\gamma_k$ must contain the geodesic $\ell$.  

 Fix the metric of constant curvature $-1$ on $N_1$, and let $U$ be an open
 neighborhood of $F_1$ in $N_1$.  For a given surface $M$, any metric
 extending  the metric surface $U$ will posses a neighborhood consisting of
 metrics, each of which has a geodesic lamination without closed ends that is
 the limit of a sequence of simple closed geodesics.
\end{proof}


\begin{thebibliography}{999}
\frenchspacing
\bibitem[Ab]{Ab}
R. Abraham,
Bumpy metrics, Global Analysis,
{\it  Proc. Sympos. Pure Math.} Vol. XIV (1968) 1--3.
\bibitem[An]{An}
D.V. Anosov, Generic properties of closed geodesics,
(Russian) {\it Izv. Akad. Nauk SSSR}, Ser. Mat. 46 (1982), no. 4, 675--709, 
896.
\bibitem[Bi]{Bi}
G.D. Birkhoff,  Dynamical Systems,  Colloquium Publ. A.M.S., Vol.9, 1927.
\bibitem[CD]{CD}
T.H. Colding and C. De Lellis, Singular limit laminations, 
Morse index, and positive scalar curvature, preprint 2002.  
\bibitem[CH1]{CH1}
T.H. Colding and N. Hingston, Metrics without Morse index bounds, 
preprint 2001.
\bibitem[CH2]{CH2}
T.H. Colding and N. Hingston, 
Morse index bounds for simple geodesics on surfaces, in preparation.
\bibitem[CM1]{CM1}
T.H. Colding and W.P. Minicozzi II,
Examples of embedded minimal tori without area bounds, {\it International
   Mathematics Research Notices}, vol. 99,
no. 20 (1999) 1097-1100.
\bibitem[CM2]{CM2}
T.H. Colding and W.P. Minicozzi II,
Embedded minimal surfaces without area bounds in $3$-manifolds,
{\it Proc. of conference on Geometry and Topology, Aarhus 1998.
Contemporary Mathematics}, vol. 258 (2000) 107-120.
\bibitem[CM3]{CM3}
T.H. Colding and W.P. Minicozzi II,
Embedded minimal surfaces in $3$-manifolds with positive scalar 
curvature, in preparation.  
\bibitem[Eb]{Eb} 
P. Eberlein, Geometry of nonpositively curved manifolds, Chicago
Lectures in Mathematics, University of Chicago Press, Chicago, Il,
1996. 
\bibitem[HaPe]{HaPe} 
J.L. Harer and R.C. Penner,   
Combinatorics of train tracks, 
Annals of Mathematics Studies, 125. 
Princeton University Press, Princeton, NJ, 1992.
\bibitem[HaNoRu]{HaNoRu}
J. Hass, P. Norbury, and J.H. Rubinstein, Minimal spheres of 
arbitrary high Morse index, preprint 2002. 
\bibitem[HiPu]{HiPu}
M.W. Hirsch and C. Pugh, Stable manifolds and
hyperbolic sets, {\it Proc. Symp. Pure Math.}, vol. 14, AMS (1970) 133-165.
\bibitem[Kl]{Kl}
W.P.A. Klingenberg,
Riemannian geometry. Second edition.
de Gruyter Studies in Mathematics, 1. Walter de Gruyter and Co., Berlin, 
1995.
\bibitem[KlT]{KlT}
W.P.A. Klingenberg and F. Takens, Generic properties of geodesic flows, 
{\it Math. Ann.} 197 (1972) 323--334.
\bibitem[MS]{MS}
D. McDuff and D. Salamon,  Introduction to Symplectic Topology.  Oxford
University Press, 1995.
\bibitem[PW]{PW} 
J. Palis and M. Welington, Geometric Theory of Dynamical Systems, 
Springer-Verlag, 1982
\bibitem[PiRu]{PiRu}
J. Pitts and J.H. Rubinstein,
Applications of minimax to minimal surfaces and
the topology
of three-manifolds,
{\it Proc. of the CMA} 12
  (1987) 137-170.
\bibitem[Sp]{Sp}
M. Spivak, A comprehensive introduction to differential geometry,
Vol. 4, chapter 8,
Second edition. Publish or Perish, Inc., Wilmington, Del., 1979.
\end{thebibliography}
\end{document}